\documentclass[12pt,]{article}
\renewcommand{\theequation}{\thesection.\arabic{equation}}

\usepackage[normal]{subfigure}
\usepackage{graphicx}
\usepackage{amssymb}
\usepackage{amscd}
\usepackage{amsmath}
\usepackage{amsfonts}
\usepackage{amsbsy}
\usepackage{epsfig}

\begin{document}
\date{}

\newenvironment{ar}{\begin{array}}{\end{array}}

\baselineskip20pt \addtocounter{page}{0}
\title{\Large\bf
Bifurcation of limit cycles by perturbing a piecewise linear Hamiltonian system with a homoclinic loop \thanks{The project supported  by  the  National Natural Science Foundation of
China (10971139) and the Slovenian Research Agency.}
 }
\author{Liang Feng$^{a, b}$\thanks {{\em E-mail
address}: liangfeng741018@163.com}, Han Maoan$^{a, c}$ \thanks{Author for correspondence. {\em E-mail address}: mahan@shnu.edu.cn.} and  Valery G. Romanovski$^{c, d}$\thanks {{\em
E-mail
address}: valery.romanovsky@uni-mb.si.} \\
{\small { $^a$The Institute of Mathematics,  Shanghai Normal University, Shanghai 200234, PR China}}\\{\small   $^b$The Institute of Mathematics,  Anhui Normal University, Wuhu 241000,
PR China}\\ {\small  $^c$Faculty of Natural Science and Mathematics,}\\
{\small University of Maribor,  SI-2000 Maribor, Slovenia}\\ {\small $^d$Center for Applied Mathematics and Theoretical Physics, }\\
{\small University of Maribor, SI-2000 Maribor, Slovenia } }

\date{ }
\maketitle

\baselineskip15pt \noindent {\bf Abstract.} In this paper, we study limit cycle bifurcations for a kind of non-smooth polynomial differential systems by perturbing a piecewise linear
Hamiltonian system with a center at the origin and a homoclinic loop around the origin. By using the first Melnikov function of piecewise near-Hamiltonian systems, we give lower bounds
of the maximal number of limit cycles in Hopf and homoclnic bifurcations, and derive an upper bound of the number of limit cycles that bifurcate from the periodic annulus between the
center and the homoclinic loop up to the first order in $\varepsilon.$ In the case  when the degree of perturbing terms is low, we obtain a precise result on the number of zeros of the
first Melnikov function.

\normalsize \baselineskip 15pt \parskip 5pt
\section{Introduction and main results}

\par~~ There are many problems in mechanics, electrical engineering
and the theory of automatic control which are described by non-smooth systems, see for instance the works of A.F.Filippov [1], Andronov et al [2], Kunze [3] and the references therein.
Recently, a good deal of work has been done to study bifurcations in non-smooth systems including Hopf, homoclinic and subharmonic bifurcations, etc.. In [4-6] Hopf bifurcation for
non-smooth systems was studied by developing new methods for computing Lyapunov constants. The Melnikov method for Hopf and homoclinic bifurcations was extended to non-smooth systems
in [7-9]. The method of averaging has also been extended to non-smooth systems in [10]. However, so far there are few papers in the literature studying homoclinic bifurcations inside
the class of piecewise polynomial differential systems. In this work, we study this problem by using the first order  Melnikov function of piecewise near-Hamiltonian systems.

\par In [7], Liu and Han considered a general form of a piecewise
near-Hamiltonian system on the plane
\begin{equation}\label{eq:100.201} \left\{\begin{array}{ll}
          \dot{x}=H_{y}+\varepsilon p(x,y,\delta),& \hbox{}\\
          \dot{y}=-H_{x}+\varepsilon q(x,y,\delta),& \hbox{}
        \end{array}
      \right.\end{equation}
where $$H(x,y)=\left\{
        \begin{array}{ll}
          H^{+}(x,y),~x\geq 0,& \hbox{}\\
            H^{-}(x,y),~x<0,& \hbox{}
        \end{array}
      \right.$$
$$p(x,y,\delta)=\left\{
        \begin{array}{ll}
          p^{+}(x,y,\delta),~x\geq 0,& \hbox{}\\
            p^{-}(x,y,\delta),~x<0,& \hbox{}
        \end{array}
      \right.
$$
$$q(x,y,\delta)=\left\{
        \begin{array}{ll}
         q^{+}(x,y,\delta),~x\geq 0,& \hbox{}\\
            q^{-}(x,y,\delta),~x<0,& \hbox{}
        \end{array}
      \right.$$
$H^{\pm},$ $p^{\pm}$ and $q^{\pm}$ are $C^{\infty}$, $\varepsilon> 0$ is small, $\delta\in D\subset R^{m}$ is a vector parameter with $D$ compact. This system has two subsystems $$
\left\{\begin{array}{ll}
          \dot{x}=H^{+}_{y}+\varepsilon p^{+}(x,y,\delta),&
          \hbox{}\\
          \dot{y}=-H^{+}_{x}+\varepsilon q^{+}(x,y,\delta),& \hbox{}
        \end{array}
      \right. \eqno(1.1a)$$
and
$$ \left\{\begin{array}{ll}
          \dot{x}=H^{-}_{y}+\varepsilon p^{-}(x,y,\delta),&
          \hbox{}\\
          \dot{y}=-H^{-}_{x}+\varepsilon q^{-}(x,y,\delta),& \hbox{}
        \end{array}
      \right. \eqno(1.1b)$$
 which are called
the right subsystem and the left subsystem,  respectively. Suppose that (\ref{eq:100.201})$|_{\varepsilon=0}$ has a family of periodic orbits around the origin and satisfies the
following two assumptions.
\par Assumption (I): There exist an interval $J=(\alpha,\beta)$, and
two points $A(h)=(0,a(h))$ and $A_{1}(h)=(0,a_{1}(h))$ such that for $h\in J$
$$H^{+}(A(h))=H^{+}(A_{1}(h))=h,~H^{-}(A(h))=H^{-}(A_{1}(h))=\widetilde{h},~a(h)\neq a_{1}(h).$$
\par Assumption (II): The subsystem
(\ref{eq:100.201}a)$|_{\varepsilon=0}$ has an orbital arc $L^{+}_{h}$ starting from $A(h)$ and ending at $A_{1}(h)$ defined by $H^{+}(x,y)=h,x\geq 0$; the subsystem
(\ref{eq:100.201}b)$|_{\varepsilon=0}$ has an orbital arc $L^{-}_{h}$ starting from $A_{1}(h)$ and ending at $A(h)$ defined by $H^{-}(x,y)=H^{-}(A_{1}(h)),x< 0$.

\par Under assumptions (I) and (II), (1.1)$|_{\varepsilon=0}$ has a family of non-smooth periodic orbits $L_{h}=L^{+}_{h}\cup L^{-}_{h},$ $h\in J.$ For definiteness, we assume the orbits $L_h$ for $h\in J$ orientate anticlockwise.
See Fig. 1.

\begin{figure}[h]
\centering {\includegraphics[bbllx=6cm,bblly=0.5cm,bburx=0.5cm, bbury=6cm,scale=0.8]{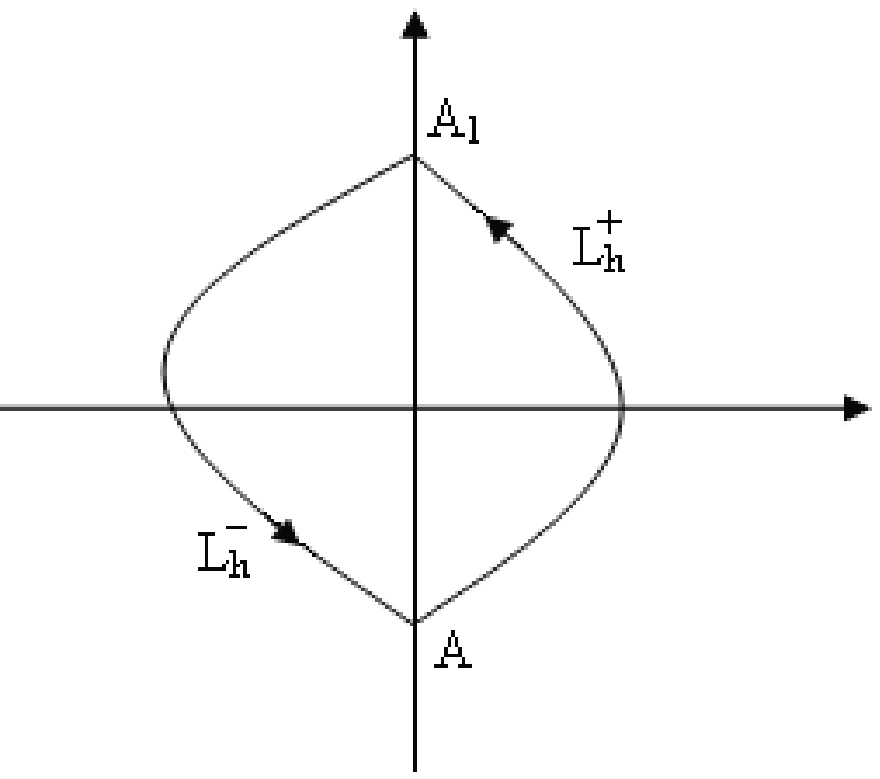}}
\begin{center}
{Fig. 1. The closed orbits of (1.1)$|_{\varepsilon=0}$.}
\end{center}
\end{figure}

 By Theorem
1.1 in [7], the first order Melnikov function of system (\ref{eq:100.201}) has the form
\begin{equation}\label{eq:100.202}
M(h,\delta)=\frac{H_{y}^{+}(A)}{H_{y}^{-}(A)}\Big[\frac{H_{y}^{-}(A_{1})}{H_{y}^{+}(A_{1})} \int_{L^{+}_{h}}q^{+}dx-p^{+}dy
 +\int_{L^{-}_{h}}q^{-}dx-p^{-}dy \Big],~h\in J.\end{equation}

Also, we know from [7] that if $M(h,\delta)$ has at most $k$ zeros in $h$ on the interval $(\alpha,\beta)$ for all $\delta\in D$, then (\ref{eq:100.201}) has at most $k$ limit cycles
bifurcated from the open annulus $\bigcup_{\alpha<h<\beta}L_h$.

Let
\begin{equation}\label{eq:100.203}
    \begin{array}{ll}
M^{+}(h,\delta)=\int_{L^{+}_{h}}q^{+}dx-p^{+}dy,~M^{-}(h,\delta)=\int_{L_{h}^{-}}q^{-}dx-p^{-}dy,&
        \hbox{}\\\\
   \displaystyle \widetilde{M}^{-}(\widetilde{h},\delta)=\int_{H^{-}(x,y)=\widetilde{h},~x\leq0}q^{-}dx-p^{-}dy,
 \end{array}\end{equation}
where $\widetilde{h}$ is given in assumption (I). Then we know
\begin{equation}\label{eq:100.1}
M^{-}(h,\delta)=\widetilde{M}^{-}(H^{-}(A_{1}(h)),\delta),~h\in J.
\end{equation}

 As in the smooth case, a very
important issue associated with (\ref{eq:100.201}) is to find the number of limit cycles and their distribution. In [7], Liu and Han consider a piecewise polynomial system of the form
\begin{equation}\label{eq:100.304}
 (\dot{x},\dot{y})=\left\{\begin{array}{ll}
          (b^{+}y+\varepsilon p^{+}_{n}(x,y,\delta),-b^{+}x+\varepsilon q^{+}_{n}(x,y,\delta)),~x\geq 0,&
          \hbox{}\\
          (b^{-}y+\varepsilon p^{-}_{n}(x,y,\delta),-b^{-}x+\varepsilon q^{-}_{n}(x,y,\delta)),~x< 0,& \hbox{}
        \end{array}
      \right.
     \end{equation}
where $b^{\pm}>0,$ $p_{n}^{\pm}$ and $q_{n}^{\pm}$ are arbitrary polynomials of degree $n$. It was proved that the maximal number of limit cycles of the above system is $n$ up to the
first order in $\varepsilon$. We know the unperturbed system of (\ref{eq:100.304}) has a global center at the origin. Then, an interesting question is how many limit cycles system
(\ref{eq:100.201}) can have if (\ref{eq:100.201})$|_{\varepsilon=0}$ is a linear piecewise system with a homoclinic loop and $p^{\pm}$ and $q^{\pm}$ are arbitrary polynomials with
degree $n$.

In this paper, we take
\begin{equation}\label{eq:100.204} H^{+}(x,y)=\frac{1}{2}\big((x-1)^2-y^2\big),~x\geq 0,\end{equation}
\begin{equation}\label{eq:100.206}H^{-}(x,y)=-\frac{1}{2}(x^2+y^2),~x<0,\end{equation}
and suppose
\begin{equation}\label{eq:100.001} \begin{array}{ll}p(x,y)=\left\{
        \begin{array}{ll}
          p^{+}(x,y)=\sum\limits_{i+j=0}^{n}a^{+}_{ij}x^{i}y^{j}, ~x\geq 0,&
          \hbox{}\\
           p^{-}(x,y)=\sum\limits_{i+j=0}^{n}a^{-}_{ij}x^{i}y^{j},~x< 0,& \hbox{}
        \end{array}
      \right.
      & \hbox{}\\
q(x,y)=\left\{
        \begin{array}{ll}
          q^{+}(x,y)=\sum\limits_{i+j=0}^{n}b^{+}_{ij}x^{i}y^{j}, ~x\geq 0,&
          \hbox{}\\
           q^{-}(x,y)=\sum\limits_{i+j=0}^{n}b^{-}_{ij}x^{i}y^{j},~x< 0.& \hbox{}
        \end{array}
      \right. \end{array}\end{equation}
Under (\ref{eq:100.204}) and (\ref{eq:100.206}), the system
 (\ref{eq:100.201}) becomes
\begin{equation}\label{eq:100.2} \left\{\begin{array}{ll}
          \dot{x}=-y+\varepsilon p^{+}(x,y),& \hbox{}\\
          \dot{y}=1-x+\varepsilon q^{+}(x,y),& \hbox{}
        \end{array} x\geq 0,
      \right.~~~~~
 \left\{\begin{array}{ll}
          \dot{x}=-y+\varepsilon p^{-}(x,y),& \hbox{}\\
          \dot{y}=x+\varepsilon q^{-}(x,y),& \hbox{}
        \end{array} x< 0.
      \right.
      \end{equation}
For system  (\ref{eq:100.2})$|_{\varepsilon=0}$, there exist a family of periodic orbits as follows
$$L_{h}=\{(x,y)|H^{+}(x,y)=\frac{h}{2},~x\geq 0\} \cup \{(x,y)|H^{-}(x,y)=\frac{\widetilde{h}}{2},~x< 0\},$$
with $\widetilde{h}=h-1,$ $0<h<1.$ For the sake of convenience, here we use $h/2$ instead of $h$. If $h\rightarrow 1,$ $L_{h}$ approaches the origin which is an elementary center of
parabolic-focus type(see [4] for the definition).
 And if $h\rightarrow 0,$ $L_{h}\rightarrow L_{0},$ where $L_{0}$ is a compound homoclinic loop with
 a saddle $S_{1}(1,0)$. See Fig. 2.
 \vspace{2.5cm}

\begin{figure}[h]
\centering {\includegraphics[bbllx=5cm,bblly=1cm,bburx=1cm,bbury=6cm,scale=0.8]{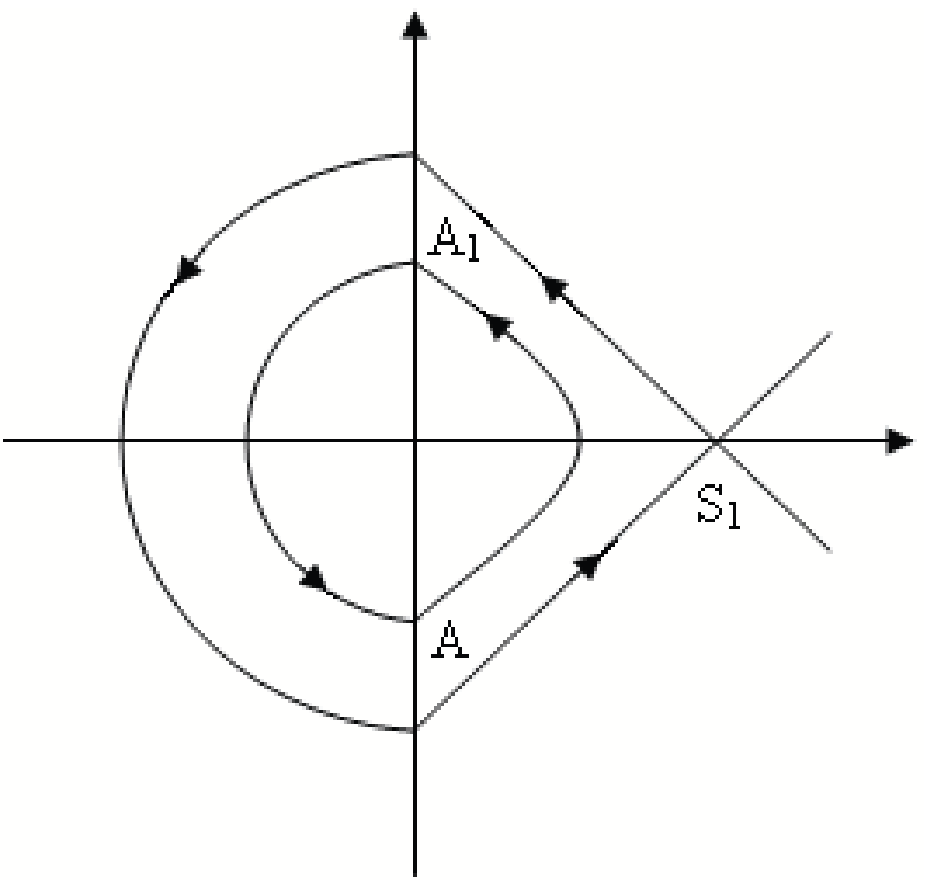}}
\begin{center}
{Fig. 2. Phase portrait of system (1.9)$|_{\varepsilon=0}$.}
\end{center}
\end{figure}

  Note that $H^{+}_{y}(0,y)\equiv H^{-}_{y}(0,y)$ for $-1<y<1$.
Then by (\ref{eq:100.202}), we have the first order Melnikov function of system
 (\ref{eq:100.2}) satisfying
 \begin{equation}\label{eq:10.3} M(\frac{h}{2})=\int_{\widehat{AA_{1}}}q^{+}dx-p^{+}dy+
 \int_{\widehat{A_{1}A}}q^{-}dx-p^{-}dy\equiv\bar M(h),\end{equation}
where $0<h<1$, and
$$\begin{array}{c}A=(0,-\sqrt{1-h}), \ A_{1}=(0,\sqrt{1-h}), \\  \widehat{AA_{1}}=\{(x,y)|H^{+}(x,y)=\frac{h}{2},~x\geq 0\},\\
 \widehat{A_{1}A}=\{(x,y)|H^{-}(x,y)=\frac{h-1}{2},~x< 0\}.\end{array}$$

Let $Z(n)$ denote the maximal number of zeros of the non-zero function $\bar M(h)$ on the open interval $(0,1)$ for all possible $p$ and $q$ satisfying \eqref{eq:100.001}, which is the
maximal number of limit cycles of (\ref{eq:100.2}) bifurcated from the periodic annulus $\bigcup_{0<h<1}L_h$ for all possible $p$ and $q$ satisfying \eqref{eq:100.001} when $\bar M(h)$
is not identically zero. Let $N_{Hopf}(n)$ and $N_{homoc}(n)$ denote, respectively,
 the maximal number of limit cycles bifurcated in Hopf bifurcation near the origin and in homoclinic bifurcation near $L_0$ for all
possible $p$ and $q$ satisfying \eqref{eq:100.001}. Then our main results can be stated as follows.

{\bf Theorem 1.1.} \ {\it For any $n\geq 1$ we have

$(1)$ $N_{Hopf}(n)\geq n+[\frac{n+1}{2}]$.

$(2)$  $N_{homoc}(n)\geq n+[\frac{n+1}{2}]$.}

{\bf Theorem 1.2.} \ {\it  For $n=1,2,3,4,$ we have $Z(n)=n+[\frac{n+1}{2}]$.}

{\bf Theorem 1.3.} \ {\it For any $n\geq 5$ we have $n+[\frac{n+1}{2}]\leq Z(n)\leq 2n+[\frac{n+1}{2}]$.}

\renewcommand\theequation{\arabic{section}.\arabic{equation}}
\makeatletter \@addtoreset{equation}{section} \makeatother
\section{Preliminary}

\par ~~~In this section we  give an expression  of the first order Melnikov function $\bar M(h)$ for $0<h<1$, together with the  two expansions of $\bar M(h)$
near the origin and the homoclinic loop $L_{0}$,  respectively. By (\ref{eq:100.203}) and (\ref{eq:100.1}) we have
\begin{equation}\label{eq:100.308}M^{+}(\frac{h}{2})=\int_{\widehat{AA_{1}}}q^{+}dx-p^{+}dy,
~~M^{-}(\frac{h}{2})=\int_{\widehat{A_{1}A}}q^{-}dx-p^{-}dy=\widetilde{M}^{-}(\frac{h-1}{2}).\end{equation} Then, applying  Green's formula twice we obtain

 $$   \begin{array}{ll}
\displaystyle M^{+}(\frac{h}{2})=-\int\int_{int(\widehat{AA_{1}}\bigcup \overrightarrow{A_{1}A})}(p^{+}_{x}+q^{+}_{y})dxdy
         -\int_{\overrightarrow{AA_{1}}}p^{+}(0,y)dy &
        \hbox{}\\\\
    ~~~~~~~~~~\displaystyle=-\oint_{\widehat{AA_{1}}\bigcup \overrightarrow{A_{1}A}}\overline{p}(x,y)dy-
            \int_{\overrightarrow{AA_{1}}}p^{+}(0,y)dy,&
        \hbox{}\\\\ \end{array}$$
where\begin{equation}\label{eq:22}\overline{p}(x,y) =p^{+}(x,y)-p^{+}(0,y)+\int_{0}^{x}q_{y}^{+}(u,y)du,\end{equation} satisfying
$$\overline{p}(0,y)=0,~~\overline{p}_{x}=p_{x}^{+}+q_{y}^{+}.$$
Since $\widehat{AA_1}$ can be represented as $x=1-\sqrt{h+y^{2}}\equiv \psi(y,h)$, we have
 \begin{equation}\label{eq:100.10}
 \begin{array}{rl}
   \displaystyle M^{+}(\frac{h}{2})&=\displaystyle -\Big[\int^{\sqrt{1-h}}_{-\sqrt{1-h}}\overline{p}(\psi(y,h),y)dy
    +\int^{\sqrt{1-h}}_{-\sqrt{1-h}}p^{+}(0,y)dy\Big]\\\\
       &\equiv -(I_{0}(h)+I_{1}(h)). \end{array}\end{equation}
Note that $p^{+}(0,y)=\sum\limits_{j=0}^{n}a^{+}_{0j}y^{j}.$ It follows that
\begin{equation}\label{eq:100.11}
    \begin{array}{ll}
\displaystyle I_{1}(h)=\int^{\sqrt{1-h}}_{-\sqrt{1-h}}p^{+}(0,y)dy&
        \hbox{}\\\\
    ~~~~~~~\displaystyle
                 =\sum\limits_{j=0}^{n}\frac{a^{+}_{0j}}{j+1}(1-(-1)^{j+1})(1-h)^{\frac{j+1}{2}}&
        \hbox{}\\\\
    ~~~~~~~\displaystyle
                 =\sqrt{1-h}\sum\limits_{k=0}^{[\frac{n}{2}]}\frac{2a^{+}_{0,2k}}{2k+1}(1-h)^{k}.
 \end{array}\end{equation}
By (\ref{eq:22}), we have
\begin{equation}\label{eq:100.101}
    \begin{array}{ll}
\displaystyle\overline{p}(x,y)=\sum\limits_{i+j=0}^{n}a^{+}_{ij}x^{i}y^{j}-\sum\limits_{j=0}^{n}a^{+}_{0j}y^{j}
                           +\sum\limits_{i+j=0}^{n}\frac{jb^{+}_{ij}}{i+1}x^{i+1}y^{j-1}&
        \hbox{}\\\\
    ~~~~~~~~~\displaystyle\equiv x\sum\limits_{i+j=0}^{n-1}p ^{+}_{ij}x^{i}y^{j},
 \end{array}\end{equation}
where
$$p ^{+}_{ij}=a ^{+}_{i+1,j}+\frac{j+1}{i+1}b ^{+}_{i,j+1}.$$ In particular, for $n=2$ we have
$$p^{+}_{00}=a^{+}_{10}+b^{+}_{01},~p^{+}_{10}=a^{+}_{20}+\frac{1}{2}b^{+}_{11},
~p^{+}_{01}=a^{+}_{11}+2b^{+}_{02}.$$
 The definition of $I_{0}(h)$ and (\ref{eq:100.101}) yield

$$
    \begin{array}{ll}
   \displaystyle I_{0}(h)=\int^{\sqrt{1-h}}_{-\sqrt{1-h}}\overline{p}(\psi(y,h),y)dy&
        \hbox{}\\\\
    ~~~~~~~\displaystyle=\sum\limits_{i+j=0}^{n-1}p^{+}_{ij}
                             \int^{\sqrt{1-h}}_{-\sqrt{1-h}}(1-\sqrt{h+y^{2}})^{i+1}y^{j}dy&
        \hbox{}\\\\\end{array}$$
\begin{equation}\label{eq:100.12}
    ~~~~~~~~~~~\displaystyle=2\sum \limits_{i+2k =0}^{n-1} p^{+}_{i,2k}
                            \int^{\sqrt{1-h}}_{0}(1-\sqrt{h+y^{2}})^{i+1}y^{2k}dy.
                       \end{equation}
For $0\leq k\leq [\frac{n-1}{2}]$ and $0\leq r\leq n$  let
\begin{equation}\label{eq:217}{I}_{rk}(h)=\int^{\sqrt{1-h}}_{0}(h+y^{2})^{\frac{r}{2}}y^{2k}dy.\end{equation} Then for $0\leq k\leq [\frac{n-1}{2}]$ and $0\leq i+2k\leq n-1$
$$\int^{\sqrt{1-h}}_{0}(1-\sqrt{h+y^{2}})^{i+1}y^{2k}dy
=\sum \limits_{r=0}^{i+1}C_{i+1}^{r}(-1)^r{I}_{rk}(h).$$ If $r=2l,~l\geq 0,$ we see
\begin{equation}\label{eq:100.13}
I_{2l,k}(h)=\int^{\sqrt{1-h}}_{0}(h+y^{2})^{l}y^{2k}dy=\sqrt{1-h}~\varphi_{l+k}(h),\end{equation}
 where
$\varphi_{l+k}(h)$ is a polynomial of  degree $l+k.$ Hence $I_{2l,k}(h)$ is $C^{\omega}$ at $h=0.$ If $r=2l+1,~l\geq 0,$ then
\begin{equation}\label{eq:100.5}I_{2l+1,k}(h)=
\int^{\sqrt{1-h}}_{0}(h+y^{2})^{l+\frac{1}{2}}y^{2k}dy.\end{equation} Using the formula
$$\int y^{2k}(h+y^{2})^{l+\frac{1}{2}}dy=\frac{y^{2k+1}(h+y^{2})^{l+\frac{1}{2}}}{2(l+k+1)}
+\frac{(2l+1)h}{2(l+k+1)}\int y^{2k}(h+y^{2})^{l-\frac{1}{2}}dy,$$ we have from (\ref{eq:100.5}) that
$$I_{2l+1,k}(h)=\frac{(1-h)^{k+\frac{1}{2}}}{2(l+k+1)}
+\frac{(2l+1)h}{2(l+k+1)}I_{2l-1,k}(h),~l\geq 1, k\geq 0.$$ It follows that

\begin{equation}\label{eq:100.6}I_{2l+1,k}(h)=(1-h)^{k+\frac{1}{2}}
\widetilde{\varphi}_{lk}(h)+\alpha_{lk}h^{l}I_{1k}(h),~l\geq 1, k\geq 0,\end{equation} where$$\alpha_{lk}=\frac{(2l+1)(2l-1)\cdots 3}{2^{l}(l+k+1)(l+k)\cdots(k+2)},$$ and
$$\widetilde{\varphi}_{lk}(h)=\frac{1}{2(l+k+1)}+\frac{(2l+1)h}{2(l+k+1)(l+k)}+\cdots
+\frac{(2l+1)(2l-1)\cdots5\cdot h^{l-1}}{2^{l}(l+k+1)(l+k)\cdots(k+2)}, $$ which is a polynomial of  degree $l-1$ in $h$. For $k\geq0,$ let
$$\varphi^{*}_{lk}(h)=\left\{
        \begin{array}{ll}
          \widetilde{\varphi}_{lk}(h), ~l\geq 1,& \hbox{}\\
           0,~~~~~~~l=0,& \hbox{}
        \end{array}
      \right.~~~~
      \alpha^{*}_{lk}=\left\{
        \begin{array}{ll}
          \alpha_{lk}, ~l\geq 1,& \hbox{}\\
           1,~~~l=0.& \hbox{}
        \end{array}
      \right.$$
By (\ref{eq:100.6}),
\begin{equation}\label{eq:100.601}I_{2l+1,k}(h)=(1-h)^{k+\frac{1}{2}}
\varphi^{*}_{lk}(h)+\alpha^{*}_{lk}h^{l}I_{1k}(h),~l\geq 0, k\geq 0.\end{equation} Further, by using the formula
$$\int y^{2k}(h+y^{2})^{\frac{1}{2}}dy=\frac{y^{2k-1}(h+y^{2})^{\frac{3}{2}}}{2(k+1)}-\frac{(2k-1)h}{2(k+1)}
\int y^{2(k-1)}(h+y^{2})^{\frac{1}{2}}dy,$$ we have
$$I_{1k}(h)=\frac{(1-h)^{k-\frac{1}{2}}}{2(k+1)}
-\frac{h(2k-1)}{2(k+1)}I_{1,k-1}(h),~k\geq 1.$$ Thus
\begin{equation}\label{eq:100.7}I_{1k}(h)=\sqrt{1-h}\psi_{k-1}(h)
+\beta_{k}h^{k}I_{10}(h),~k\geq 1,\end{equation} where$$\beta_{k}=(-1)^{k}\frac{(2k-1)!!}{2^{k}(k+1)!},$$ and
$$\psi_{k-1}(h)=\frac{(1-h)^{k-1}}{2(k+1)}-\frac{(2k-1)h(1-h)^{k-2}}{4(k+1)k}+\cdots+(-1)^{k-1}\frac{(2k-1)(2k-3)
\cdots 3 \cdot h^{k-1}} {2^{k}(k+1)k(k-1)\cdots2},$$ which is a polynomial in $h$ of  degree $k-1$.  Let
$$\psi^{*}_{k-1}(h)=\left\{
        \begin{array}{ll}
         \psi_{k-1}(h), ~k\geq 1,& \hbox{}\\
           0,~~~~~~~k=0,& \hbox{}
        \end{array}
      \right.~~~~
      \beta^{*}_{k}=\left\{
        \begin{array}{ll}
          \beta_{k}, ~k\geq 1,& \hbox{}\\
           1,~~~k=0.& \hbox{}
        \end{array}
      \right.$$
By (\ref{eq:100.7}),
\begin{equation}\label{eq:100.701}I_{1k}(h)=\sqrt{1-h}\psi^{*}_{k-1}(h)
+\beta^{*}_{k}h^{k}I_{10}(h),~k\geq 0.\end{equation} Moreover, by (\ref{eq:100.5})

\begin{equation}\label{eq:100.16}
    \begin{array}{ll}
  \displaystyle I_{10}(h)=\int^{\sqrt{1-h}}_{0}(h+y^{2})^{\frac{1}{2}}dy
=\frac{1}{2}[\sqrt{1-h}+h ~{\rm ln}(1+\sqrt{1-h})-\frac{1}{2}h~{\rm ln} h]
\equiv-\frac{1}{4}h{\rm ln}h+\gamma(h),
                       \end{array}\end{equation}
where $\gamma (h)$ is $C^{\omega}$ at $h=0.$ We observe that \begin{equation}\label{eq:224} I_{10}(h)=h\varphi_0\Big(\sqrt{\frac{1-h}{h}}\Big),\
\varphi_0(u)=\int_0^{u}\sqrt{1+x^2}dx.\end{equation} The function $ \varphi_0(u)$ is analytic on $R$ and odd in $u$.

 From (\ref{eq:100.601}), (\ref{eq:100.701}) and (\ref{eq:100.16}), we can obtain
\begin{equation}\label{eq:100.8}
    \begin{array}{rl}
  I_{2l+1,k}(h)&=\displaystyle(1-h)^{k+\frac{1}{2}}\varphi^{*}_{lk}(h)
+\sqrt{1-h}\alpha^{*}_{lk}\psi^{*}_{k-1}(h)h^{l}+\alpha^{*}_{lk}\beta^{*}_{k}h^{l+k}I_{10}(h)\\\\
&=\displaystyle\sqrt{1-h}\bar \varphi_{l+k-1}(h)+\bar \alpha_{lk}h^{l+k}I_{10}(h)\\\\
     &=\displaystyle\widetilde{\gamma}(h)-\frac{1}{4}\bar \alpha_{lk}h^{l+k+1}{\rm ln}h,
                       \end{array}\end{equation}
where $\bar \varphi_{l+k-1}(h)$ is a polynomial of degree $l+k-1$, $\bar \alpha_{lk}$ are non-zero constants,
  $\widetilde{\gamma} (h)$ is $C^{\omega}$ at $h=0,$
$0\leq k \leq[\frac{n-1}{2}],$  $0\leq 2l+1\leq n$.

Substituting  (\ref{eq:100.13}) and (\ref{eq:100.8}) into (\ref{eq:100.12}), and noting that $k+[\frac{i+1}{2}]\leq [\frac{n}{2}]$ and  $k+[\frac{i}{2}]\leq [\frac{n-1}{2}]$ since
$i+2k\leq n-1$, we have

\begin{equation}\label{eq:100.9}
    \begin{array}{ll}
\displaystyle I_{0}(h)&=2\sum \limits_{i+2k =0}^{n-1} p^{+}_{i,2k}\Big(\sum \limits_{r=0,~r~ {\rm even}}^{i+1}C_{i+1}^{r}I_{rk}-\sum \limits_{r=1,~r{\rm~
odd}}^{i+1}C_{i+1}^{r}I_{rk}\Big)\\\\
     &=\displaystyle2\sum \limits_{i+2k =0}^{n-1} p^{+}_{i,2k}
    \Big[\sqrt{1-h}\nu_{[\frac{i+1}{2}]+k}(h)+ \nu_{[\frac{i}{2}]}(h)h^kI_{10})\Big]\\\\
    &=\sqrt{1-h} \nu_{[\frac{n}{2}]}(h)+ \nu_{[\frac{n-1}{2}]}(h)I_{10}  \\\\
    &\equiv\displaystyle \phi_{1}(h){\rm ln}h+\phi_{2}(h),\end{array}\end{equation}
where $\nu_j(h)$ denotes a polynomial of $h$ of  degree $j$, $\phi_{2}(h)\in C^{\omega}$ at $h=0$, $\phi_{1}(h)$ is a polynomial in $h$ of  degree $[\frac{n+1}{2}] $ and
$\phi_{1}(0)=0$.

Now we consider $M^{-}(\frac{h}{2})$ in (\ref{eq:100.308}).  Set $x=\sqrt{-\widetilde{h}}{\rm ~cos}\theta,~y=\sqrt{-\widetilde{h}}{~\rm sin}\theta$ for
$-1<\widetilde{h}<0,~\frac{\pi}{2}\leq\theta\leq\frac{3\pi}{2},$  and let
$$\bar I_{ij}=\int^{\frac{3\pi}{2}}_{\frac{\pi}{2}}{\rm cos}^{i}\theta~{\rm
sin}^{j}\theta d\theta
        =(-1)^{i+j}\int^{\frac{\pi}{2}}_{\frac{-\pi}{2}}{\rm cos}^{i}\theta~{\rm
        sin}^{j}\theta d\theta
        =(-1)^{i+j}\Gamma_{ij},$$
where
$$\Gamma_{ij}=\left\{
        \begin{array}{ll}
          0,~j ~{\rm odd},& \hbox{}\\
            {\rm a ~positive ~constant},~j~{\rm even}.& \hbox{}
        \end{array}
      \right.$$
It follows from (\ref{eq:100.203}) that
\begin{equation}\label{eq:100.41}
    \begin{array}{ll}
    ~~~\displaystyle\widetilde{M}^{-}(\frac{\widetilde{h}}{2})=\int_{H^{-}(x,y)=\frac{\widetilde{h}}{2}}q^{-}dx-p^{-}dy&
        \hbox{}\\\\
    ~~~~~~~~~~~~~\displaystyle=-\sum\limits_{i+j=0}^{n}(-\widetilde{h})^{\frac{i+j+1}{2}}
                             \int^{\frac{3\pi}{2}}_{\frac{\pi}{2}}(b^{-}_{ij}{\rm
                             cos}^{i}\theta~
                             {\rm sin}^{j+1}\theta+a^{-}_{ij}{\rm
                             cos}^{i+1}\theta~
                             {\rm sin}^{j}\theta)d\theta &
        \hbox{}\\\\
    ~~~~~~~~~~~~~\displaystyle=-\sum\limits_{i+j=0}^{n}(-\widetilde{h})^{\frac{i+j+1}{2}}
                             (b^{-}_{ij}\bar I_{i,j+1}+a^{-}_{ij}\bar I_{i+1,j})&
        \hbox{}\\\\
    ~~~~~~~~~~~~~\displaystyle=-\sum\limits_{l=0}^{n}(-\widetilde{h})^{\frac{l+1}{2}}\sum\limits_{j=0}^{l}
                             (b^{-}_{l-j,j}\bar I_{l-j,j+1}+a^{-}_{l-j,j}\bar I_{l-j+1,j})&
        \hbox{}\\\\
    ~~~~~~~~~~~~~\displaystyle=-\sum\limits_{l=0}^{n}e_{l}(-\widetilde{h})^{\frac{l+1}{2}},
                       \end{array}\end{equation}
where
 \begin{equation}\label{eq:100.420} \begin{array}{ll}
 \displaystyle e_{l}=\sum\limits_{j=0}^{l}(b^{-}_{l-j,j}\bar I_{l-j,j+1}+a^{-}_{l-j,j}\bar I_{l-j+1,j})&
        \hbox{}\\\\
   ~~\displaystyle=(-1)^{l+1}\sum\limits_{j=0}^{l}(b^{-}_{l-j,j}\Gamma_{l-j,j+1}+a^{-}_{l-j,j}\Gamma_{l-j+1,j})&
        \hbox{}\\\\
    ~~\displaystyle=(-1)^{l+1}\big[\sum\limits_{j=0,~j~{\rm
even}}^{l} a^{-}_{l-j,j}\Gamma_{l-j+1,j}+\sum\limits_{j=1,~j~{\rm odd}}^{l}b^{-}_{l-j,j}\Gamma_{l-j,j+1}\big] , ~0\leq l\leq n.
\end{array}\end{equation}
By (\ref{eq:100.308}) and (\ref{eq:100.41}),
\begin{equation}\label{eq:100.42}
M^{-}(\frac{h}{2})=-\sqrt{1-h}\sum\limits_{l=0}^{n}e_{l}(1-h)^{\frac{l}{2}},~~~0<h<1,
\end{equation}
where for $n=2$
$$e_{0}=-2a^{-}_{00},~e_{1}=\frac{\pi}{2}(a^{-}_{10}+b^{-}_{01}),~e_{2}=-\frac{2}{3}(2a^{-}_{20}-b^{-}_{11}+a^{-}_{02}).$$

Hence by \eqref{eq:10.3}, \eqref{eq:100.308}, \eqref{eq:100.10},  \eqref{eq:100.11},  \eqref{eq:100.9} and  \eqref{eq:100.42} we obtain the following lemma.

{\bf Lemma 2.1.}\ {\it For system \eqref{eq:100.2}, the first order Melnikov function has the following form

$$\bar M(h)=\sqrt{1-h} f_n(\sqrt{1-h})+g_
{[\frac{n-1}{2}]}(1-h)I_{10}(h),$$ where $ f_n$ and $g_ {[\frac{n-1}{2}]}$ are polynomials of degrees  $n$ and $[\frac{n-1}{2}]$ respectively.}

\par In the following we study the expansions of $\bar M(h)$ near $h=0,1$. By (\ref{eq:100.16}) and Lemma 2.1, the lemma below holds immediately.
\par {\bf Lemma 2.2.}\ {\it For system \eqref{eq:100.2}
 the first order  Melnikov function $\bar{M}(h)$ has the following expansion near the homoclinic loop $L_{0}$
$$\bar{M}(h)=(\sum\limits_{i=1}^{[\frac{n+1}{2}]}b^{*}_{i}h^{i}){~\rm ln}h
                  +\sum\limits_{j\geq 0}b_{j}h^{j},~~0<h\ll 1,$$
where for $n=2$
$$  \begin{array}{ll}
   b_{1}^{*}=-\frac{1}{2}(\,{\it a}^{+}_{{10}}+\,{\it b}^{+}_{{01}}+2{\it a}^{+}_{{20}}+\,{\it b}^{+}_{{11}}),&
        \hbox{}\\\\
   b_{0}=-{\it a}^{+}_{{10}}-{\it b}^{+}_{{01}}-\frac{2}{3}\,{\it a}^{+}_{{20}}-\frac{1}{3}\,{\it b}^{+}_{{11}}
         -2\,{\it a}^{+}_{{00}}-\frac{2}{3}\,{\it a}^{+}_{{02}}+2a^{-}_{00}-\frac{\pi}{2}a^{-}_{10}-\frac{\pi}{2}b^{-}_{01}
         +\frac{4}{3}a^{-}_{20}-\frac{3}{2}b^{-}_{11}+\frac{2}{3}a^{-}_{02},&
        \hbox{}\\\\
  b_{1}=(\frac{1}{2}+{\rm ln}2)\,{\it a}^{+}_{{10}}+(\frac{1}{2}+{\rm ln}2){\it b}^{+}_{{01}}
         +(2{\rm ln}2-1){\it a}^{+}_{{20}}+({\rm ln}2-\frac{1}{2})\,{\it
         b}^{+}_{{11}} +{\it a}^{+}_{{00}}+{\it a}^{+}_{{02}}-a^{-}_{00}&
        \hbox{}\\\\
        ~~~~~~+\frac{\pi}{2}a^{-}_{10}+
         \frac{\pi}{2}b^{-}_{01}-2a^{-}_{20}+b^{-}_{11}-a^{-}_{02},&
        \hbox{}\\\\
    b_{2}=\frac{1}{4}\,{\it a}^{+}_{{20}}+\frac{1}{8}\,{\it b}^{+}_{{11}}-\frac{1}{4}\,{\it a}^{+}_{{02}}
       -\frac{1}{8}\,{\it a}^{+}_{{10}}-\frac{1}{8}\,{\it b}^{+}_{{01}}+\frac{1}{4}\,{\it a}^{+}_{{00}}
       -\frac{1}{4}a^{-}_{00}+\frac{1}{2}a^{-}_{02}-\frac{1}{4}b^{-}_{11}+\frac{1}{4}a^{-}_{02},&
        \hbox{}\\\\ ~\cdots\cdots
                       \end{array}$$
                       }
\par Noticing that in (\ref{eq:224}) $ \varphi_0(u)\in C^{w}$ on $R$ and is odd in $u$, we can write for $|u|$ small
$$\varphi_{0}(u)=\sum\limits_{i=0}^{\infty}\varsigma_{i}u^{2i+1},$$
where $\varsigma_{i}$ is a constant, $i\geq 0.$ Then, it follows from (\ref{eq:224}) for $1-h>0$ small
\begin{equation}\label{eq:231}I_{10}(h)=h\varphi_0\Big(\sqrt{\frac{1-h}{h}}\Big)
                =\sqrt{h(1-h)}\sum\limits_{i=0}^{\infty}\varsigma_{i}(\frac{1-h}{h})^{i}\equiv \sqrt{1-h}\phi_{3}(1-h),\end{equation}
where $\phi_{3}(u)\in C^{w}$ at $u=0.$ Following (\ref{eq:231}) and Lemma 2.1, we have
\par {\bf Lemma 2.3.}\ {\it For system (\ref{eq:100.2})
 the first order  Melnikov function $\bar{M}(h)$ has the following expansion near the origin
$$
\bar{M}(h)=\sqrt{1-h}\big[\sum\limits_{i=0}^{n}c_{i}(1-h)^{\frac{i}{2}}
                      +\sum\limits_{j\geq 1}c_{n+j}(1-h)^{[\frac{n}{2}]+j}\big],~~0<1-h\ll1,
$$
where for $n=2$
$$
    \begin{array}{ll}
    c_{0}=2(a^{-}_{00}-a^{+}_{00}),~c_{1}=-\frac{\pi}{2}(a^{-}_{10}+b^{-}_{01}),&
        \hbox{}\\\\
    c_{2}=-\frac{2}{3}(a^{+}_{02}+b^{+}_{01}+a^{+}_{10}-2a^{-}_{20}+b^{-}_{11}-a^{-}_{02}),&
        \hbox{}\\\\
   c_{3}=-\frac{2}{15}(2a^{+}_{20}+a^{+}_{10}+b^{+}_{11}+b^{+}_{01}),\cdots \cdots                           \end{array}$$
     }

\baselineskip 16pt \parskip 5pt
\section{ Proof of the main results}

\par {\bf Proof of Theorem 1.1.} {\bf (1) } For simplicity, we let $p^{\pm}$ and $q^{\pm}$ in (\ref{eq:100.001}) satisfy
\begin{equation}\label{eq:3.301}
p^{\pm}(x,y)=\sum\limits_{i= 0}^{n}a^{\pm}_{i}x^{i},~q^{\pm}(x,y)\equiv 0,
\end{equation}
that is,  $p^{\pm}(x,y)$ are univariate polynomials of the variable  $x$. By (\ref{eq:100.11}) and (\ref{eq:100.101}), we have
$$I_{1}(h)=2a^{+}_{0}\sqrt{1-h},~~\overline{p}(x,y)=\sum\limits_{i=
1}^{n}a_{i}^{+}x^{i},$$ From (\ref{eq:100.204}), the equation $H^{+}(x,y)=\frac{h}{2}$ is equivalent to
\begin{equation}\label{eq:100.21} 2x-x^{2}=\omega,\end{equation}
where $\omega =1-h-y^{2}.$ Consider equation (\ref{eq:100.21}) in $x$ near $x=0.$
 There is a unique $C^{\infty}$ solution
\begin{equation}\label{eq:100.22} x=\sigma(\omega)=\frac{1}{2}\omega+O(\omega^{2})\equiv\sum\limits_{j\geq 1}g_{j}\omega^{j},\end{equation}
where $$g_{1}=\frac{1}{2},~g_{2}=\frac{1}{8},~g_{3}=\frac{1}{16},~g_{4}=\frac{5}{128},\cdots $$ Then,
 (\ref{eq:100.22})  implies  that  the arc $\widehat{AA_1}$ can be
represented by $x=\sigma(\omega)$ near the origin. That is, in this case the function $\psi(y,h)$ in (\ref{eq:100.10}) can be taken as $x=\sigma(\omega)$ with $\omega
=1-h-y^{2},~1-h>0$ small. Note that
$$\overline{p}(\sigma(\omega),y)=\sum\limits_{i=
1}^{n}a_{i}^{+}(\sigma(\omega))^{i}=\sum\limits_{j\geq1}\widetilde{p}_{j}\omega^{j},$$ where $\widetilde{p}_{1}=\frac{1}{2}a_{1}^{+},$
$\widetilde{p}_{2}=\frac{1}{4}a_{2}^{+}+\frac{1}{8}a_{1}^{+},$ $\cdots,$ $\widetilde{p}_{n}=\frac{a_{n}^{+}}{2^{n}}+L(a^{+}_{1},\cdots,a^{+}_{n-1}),\cdots$,
 throughout this paper $L(\cdot)$ denotes  a
linear combination. Then by (\ref{eq:100.12}), we have for $0<1-h\ll 1$
$$ \begin{array}{ll}
    \displaystyle I_{0}(h)=\sum\limits_{j\geq1}\int^{\sqrt{1-h}}_{-\sqrt{1-h}}\widetilde{p}_{j}\omega^{j}dy&
        \hbox{}\\\\
    ~~~~~~~\displaystyle=2\sum\limits_{j\geq1}\widetilde{p}_{j}\int_{0}^{\sqrt{1-h}}(1-h-y^{2})^{j}dy&
        \hbox{}\\\\
    ~~~~~~~\displaystyle=\sqrt{1-h}\sum\limits_{j\geq1}\widetilde{p}_{j}B_{j}(1-h)^{j},
                       \end{array}$$
where $B_{j}=2\int_{0}^{1}(1-u^{2})^{j}du$ is a positive constant for $j\geq 1.$ Hence, by (\ref{eq:100.10})
\begin{equation}\label{eq:3.901} \begin{array}{ll}
    \displaystyle M^{+}(\frac{h}{2})=-\big(I_{0}(h)+I_{1}(h)\big)&
        \hbox{}\\\\
    ~~~~~~~~~~\displaystyle=-\big(\sqrt{1-h}\sum\limits_{j\geq1}\widetilde{p}_{j}B_{j}(1-h)^{j}+2a^{+}_{0}\sqrt{1-h}\big)&
        \hbox{}\\\\
    ~~~~~~~~~~\displaystyle\equiv\sqrt{1-h} \sum\limits_{j\geq 0} c^{*}_{j}(1-h)^{j},~{\rm for} ~0<1-h\ll 1,
                 \end{array}\end{equation}
where
 $c^{*}_{j}=\left\{\begin{array}{ll}
          -2a_{0}^{+},~j=0,& \hbox{}\\
          -\widetilde{p}_{j}B_{j},~~j\geq 1. & \hbox{}
        \end{array}
      \right.$
Under (\ref{eq:3.301}), we also have from (\ref{eq:100.420}) and (\ref{eq:100.42})
\begin{equation}\label{eq:3.902} M^{-}(\frac{h}{2})=
-\sqrt{1-h} \sum\limits_{l=0}^{n} e_{l}(1-h)^{\frac{l}{2}},~ ~0<h< 1,\end{equation} where $e_{l}=(-1)^{l+1}a_{l}^{-}\Gamma_{l+1,0},$ $\Gamma_{l+1,0}$ is a positive constant, $0\leq
l\leq n.$ Thus, by (\ref{eq:10.3}), (\ref{eq:3.901}) and (\ref{eq:3.902})
\begin{equation}\label{eq:3.903} \bar{M}(h)=
\sqrt{1-h} \big[\sum\limits_{i=0}^{n} c_{i}(1-h)^{\frac{i}{2}}+\sum\limits_{j\geq1} c_{n+j}(1-h)^{[\frac{n}{2}]+j}\big],~ ~0<1-h\ll 1,\end{equation} with
$$  \begin{array}{ll}
    \displaystyle
   c_{0}=c^{*}_{0}-e_{0}=-2a_{0}^{+}+\Gamma_{10}a_{0}^{-},~~~~
    \displaystyle
    c_{1}=-e_{1}=-\Gamma_{20}a_{1}^{-},&
        \hbox{}\\\\
    \displaystyle
    c_{2}=c^{*}_{1}-e_{2}=-\frac{B_{1}}{2}a_{1}^{+}+\Gamma_{30}a_{2}^{-},~~~~
    \displaystyle
    c_{3}=-e_{3}=-\Gamma_{40}a_{3}^{-},~~    \cdots\cdots &
    \hbox{}\\\\
\displaystyle c_{n}=\left\{
        \begin{array}{ll}
          c^{*}_{\frac{n}{2}}-e_{n}=-\frac{B_{\frac{n}{2}}}{2^{\frac{n}{2}}}a^{+}_{\frac{n}{2}}+(-1)^{n}\Gamma_{n+1,0}a_{n}^{-}
           +L(a_{1}^{+},\cdots,a^{+}_{\frac{n}{2}-1}),~n ~{\rm even},&
           \hbox{}\\\\
            -e_{n}=(-1)^{n}\Gamma_{n+1,0}a_{n}^{-},~n~{\rm odd},& \hbox{}
        \end{array}
      \right.&
        \hbox{}\\\\   \displaystyle c_{n+1}=c^{*}_{[\frac{n}{2}]+1}=
    -\frac{B_{[\frac{n}{2}]+1}}{2^{[\frac{n}{2}]+1}}a_{[\frac{n}{2}]+1}^{+}+L(a_{1}^{+},\cdots,a^{+}_{[\frac{n}{2}]}),&
    \hbox{}\\\\
    \cdots\cdots &
    \hbox{}\\\\
    \displaystyle c_{2n-[\frac{n}{2}]}=c^{*}_{n}=
    -\frac{B_{n}}{2^{n}}a_{n}^{+}+L(a_{1}^{+},\cdots,a^{+}_{n-1}),&
    \hbox{}\\\\
    \displaystyle \cdots\cdots
                       \end{array}~~~~~~~~~~~~~~~~~~~~~~~~~~~~~~~~~~~~~~~~~~~$$
Following the formulas above, we have {\small $$ \begin{array}{ll}
    ~~~~~~~~~~~~~~~~~~~~~~~~\displaystyle \frac{\partial(c_{0},c_{1},c_{2},\cdots,
    c_{n},c_{n+1},\cdots,c_{2n-[\frac{n}{2}]})}
{\partial(a^{+}_{0},a^{+}_{1},\cdots,a^{+}_{[\frac{n}{2}]},a^{+}_{[\frac{n}{2}]+1},\cdots,a^{+}_{n},a^{-}_{0},a^{-}_{1}, \cdots,a^{-}_{n})}= &
    \hbox{}\\\\
    \displaystyle \left(
   \begin{array}{ccccccccccccc}
     -2 & 0 & \cdots & 0 & 0 & 0 & \cdots & 0 & \Gamma_{10}& 0 &\cdots & 0\\
      0 & 0 & \cdots & 0 & 0 & 0 & \cdots & 0 & 0 & -\Gamma_{20} &\cdots & 0\\
      0 & -\frac{B_{1}}{2} & \cdots & 0 & 0 & 0 & \cdots & 0 & 0 & 0 &\cdots & 0\\
     \cdots & \cdots & \cdots & \cdots & \cdots & \cdots & \cdots & \cdots & \cdots & \cdots & \cdots & \cdots\\
     0 & * & \cdots & X^{*} & 0 & 0 & \cdots & 0 & 0 & 0 & \cdots & (-1)^{n}\Gamma_{n+1,0}\\
     0 & * & \cdots & * & -\frac{B_{[\frac{n}{2}]+1}}{2^{[\frac{n}{2}]+1}} & 0 & \cdots & 0 & 0 & 0 & \cdots & 0\\
     0 & * & \cdots & * & * & -\frac{B_{[\frac{n}{2}]+2}}{2^{[\frac{n}{2}]+2}} & \cdots & 0 & 0 & 0 & \cdots & 0\\
     \cdots & \cdots & \cdots & \cdots & \cdots & \cdots & \cdots & \cdots & \cdots & \cdots & \cdots & \cdots\\
      0 & * & \cdots & * & * & * & \cdots & -\frac{B_{n}}{2^{n}} & 0 & 0 & \cdots & 0\\
       \end{array}
 \right)
 \end{array}$$}
where $$X^{*}=\left\{
        \begin{array}{ll}
          -\frac{B_{\frac{n}{2}}}{2^{\frac{n}{2}}},~n ~{\rm even},&
           \hbox{}\\
            0,~n~{\rm odd}.& \hbox{}
        \end{array}
      \right.$$
Clearly, the rank of this matrix is $2n-[\frac{n}{2}]+1,$ which means that $c_{i},~0\leq i\leq 2n-[\frac{n}{2}]$ can be chosen as free parameters such that $0<|c_{0}|\ll
|c_{1}|\ll|c_{2}|\ll\cdots\ll|c_{n}|\ll|c_{n+1}|\ll|c_{n+2}|\ll\cdots \ll|c_{2n-[\frac{n}{2}]}|\ll 1,$ and $c_{i}c_{i+1}<0,~0\leq i\leq 2n-[\frac{n}{2}]-1.$ Hence, by (\ref{eq:3.903}),
in this case $\bar{M}(h)$ has $2n-[\frac{n}{2}]$ positive simple zeros in $(0,1)$  near   $h=1.$ Therefore, system (\ref{eq:100.2}) can have $2n-[\frac{n}{2}]$
(=$n+[\frac{n+1}{2}],~n\geq 1$) limit cycles near the origin, which means $N_{Hopf}(n)\geq n+[\frac{n+1}{2}]$.

\par {\bf (2)} Here, we suppose  $p^{\pm}$ and $q^{\pm}$ in (\ref{eq:100.001}) satisfy
\begin{equation}\label{eq:3.302}
p^{+}(x,y)=\sum\limits_{i= 0}^{n}a_{i}^{+}y^{i},~q^{+}(x,y)=\sum\limits_{i= 0}^{n}b_{i}^{+}y^{i},~ p^{-}(x,y)=\sum\limits_{i=0}^{n}a^{-}_{i}x^{i},~q^{-}(x,y)\equiv 0,
\end{equation}
where $p^{+}(x,y)$ and $q^{+}(x,y)$ are independent of $x$ and $p^{-}(x,y)$ is independent of $y$ . Then, from (\ref{eq:100.11}) to (\ref{eq:100.16}), we have
$$\overline{p}(x,y)=x\sum\limits_{i=1}^{n}b^{+}_{i}iy^{i-1}=x\sum\limits_{j=0}^{n-1}p^{+}_{0j}y^{j},$$
where
\begin{equation}\label{eq:3.601} p^{+}_{00}=b^{+}_{1},p^{+}_{01}=2b^{+}_{2},\cdots,p^{+}_{0,n-1}=nb^{+}_{n}.
\end{equation}
Moreover,
$$ I_{1}(h)=\sqrt{1-h}\sum\limits_{k=0}^{[\frac{n}{2}]}\frac{2a^{+}_{2k}}{2k+1}(1-h)^{k},$$ and
$$  \begin{array}{ll}
    \displaystyle I_{0}(h)=\sum\limits_{j=0}^{n-1}p^{+}_{0j}\int^{\sqrt{1-h}}_{-\sqrt{1-h}}(1-\sqrt{h+y^{2}})y^{j}dy&
        \hbox{}\\\\
    ~~~~~~~\displaystyle=2\sum\limits_{k=0}^{[\frac{n-1}{2}]}p^{+}_{0,2k}
    \int^{\sqrt{1-h}}_{0}(1-\sqrt{h+y^{2}})y^{2k}dy&
        \hbox{}\\\\
    ~~~~~~~\displaystyle=2\Big[\sum\limits_{k=0}^{[\frac{n-1}{2}]}p^{+}_{0,2k}
    \int^{\sqrt{1-h}}_{0}y^{2k}dy-\sum\limits_{k=0}^{[\frac{n-1}{2}]}p^{+}_{0,2k}
    \int^{\sqrt{1-h}}_{0}\sqrt{h+y^{2}}y^{2k}dy   \Big]&
        \hbox{}\\\\\end{array}$$
    $$\displaystyle=2\Big[\sqrt{1-h}\sum\limits_{k=0}^{[\frac{n-1}{2}]}\frac{p^{+}_{0,2k}}{2k+1}(1-h)^{k}
    -\sum\limits_{k=0}^{[\frac{n-1}{2}]}p^{+}_{0,2k}I_{1k}(h)\Big],~~~~~~~~$$

where $I_{1k}$ is given by (\ref{eq:217}). Hence by (\ref{eq:100.10}), (\ref{eq:100.701}), (\ref{eq:100.16}) and the above
  $$\begin{array}{ll}
\displaystyle \frac{M^{+}(\frac{h}{2})}{-2\sqrt{1-h}}&=\displaystyle\sum\limits_{k=0}^{[\frac{n}{2}]}\frac{a^{+}_{2k}}{2k+1}(1-h)^{k}
    +\sum\limits_{k=0}^{[\frac{n-1}{2}]}\frac{p^{+}_{0,2k}}{2k+1}(1-h)^{k}
    -\sum\limits_{k=0}^{[\frac{n-1}{2}]}p^{+}_{0,2k}\psi^{*}_{k-1}(h)
        \hbox{}\\\\
  & ~~\displaystyle-(1-h)^{-\frac{1}{2}}\sum\limits_{k=0}^{[\frac{n-1}{2}]}p^{+}_{0,2k}\beta^{*}_{k}h^{k}I_{10}(h)
   \end{array}$$
   \begin{equation}\label{eq:3.6} ~~~~~~~~~~~~~~~~~~~~~~~~~~\begin{array}{ll}
        &\displaystyle=\sum\limits_{k=0}^{[\frac{n}{2}]}\frac{a^{+}_{2k}}{2k+1}(1-h)^{k}
    +\sum\limits_{k=0}^{[\frac{n-1}{2}]}\frac{p^{+}_{0,2k}}{2k+1}(1-h)^{k}
    -\sum\limits_{k=0}^{[\frac{n-1}{2}]}p^{+}_{0,2k}\psi^{*}_{k-1}(h)
            \hbox{}\\\\
     &\displaystyle~~~~-\frac{1}{2}\sum\limits_{k=0}^{[\frac{n-1}{2}]}p^{+}_{0,2k}\beta^{*}_{k}h^{k}
    -\frac{1}{2}\sum\limits_{k=0}^{[\frac{n-1}{2}]}p^{+}_{0,2k}\beta^{*}_{k}h^{k+1}(1-h)^{-\frac{1}{2}}{\rm
    ln}(1+\sqrt{1-h})
            \hbox{}\\\\
     &\displaystyle~~~~+\frac{1}{4}\sum\limits_{k=0}^{[\frac{n-1}{2}]}p^{+}_{0,2k}\beta^{*}_{k}(1-h)^{-\frac{1}{2}}h^{k+1}{\rm
    ln}h
        \hbox{}\\\\
    &\displaystyle\equiv \sum\limits_{k=0}^{\infty}v_{k}h^{k}+\sum\limits_{k=0}^{\infty}v^{*}_{k}h^{k+1}{\rm
    ln}h,~~{\rm for}~0<h\ll 1,
                       \end{array}\end{equation}

where
$$  \begin{array}{ll}
    \displaystyle
    v_{0}=a_{0}^{+}+L(a_{2}^{+},\cdots,a_{2[\frac{n}{2}]}^{+},p_{00}^{+},\cdots,p_{0,2[\frac{n-1}{2}]}^{+}),&
        \hbox{}\\\\
    \displaystyle v_{1}=-\frac{a_{2}^{+}}{3}+L(a_{4}^{+},\cdots,a_{2[\frac{n}{2}]}^{+},
                                  p_{00}^{+},\cdots,p_{0,2[\frac{n-1}{2}]}^{+}),&
        \hbox{}\\\\
    \displaystyle v_{2}=\frac{a_{4}^{+}}{5}+L(a_{6}^{+},\cdots,a_{2[\frac{n}{2}]}^{+},
                                  p_{00}^{+},\cdots,p_{0,2[\frac{n-1}{2}]}^{+}), ~~\cdots\cdots &
    \hbox{}\\\\
    \displaystyle v_{[\frac{n}{2}]}=\frac{a_{2[\frac{n}{2}]}^{+}}{2[\frac{n}{2}]+1}(-1)^{[\frac{n}{2}]}
                                          +L(p_{00}^{+},\cdots,p_{0,2[\frac{n-1}{2}]}^{+}),&
    \hbox{}\\\\
    \displaystyle    v_{0}^{*}=\frac{1}{4}p_{00}^{+},~ v_{1}^{*}=\frac{1}{4}p_{02}^{+}\beta_{1}^{*}+L(p_{00}^{+}),&
        \hbox{}\\\\
    \displaystyle
    v_{2}^{*}=\frac{1}{4}p_{04}^{+}\beta_{2}^{*}+L(p_{00}^{+},p_{02}^{+}),&
    \hbox{}\\\\
    \displaystyle \cdots\cdots \end{array}
$$\begin{equation}\label{eq:3.611}
    \displaystyle v^{*}_{[\frac{n-1}{2}]}=\frac{1}{4}p_{0,2[\frac{n-1}{2}]}^{+}\beta_{[\frac{n-1}{2}]}^{*}
            +L(p_{00}^{+},p_{02}^{+},\cdots,p_{0,2[\frac{n-1}{2}]-2}^{+}).~~~~~
                       \end{equation}
Moreover, from (\ref{eq:3.902}) we have for $~0<h\ll 1$
\begin{equation}\label{eq:4.9}
\frac{M^{-}(\frac{h}{2})}{-2\sqrt{1-h}}=\frac{1}{2}\sum\limits_{l=0}^{n} e_{l}(1-h)^{\frac{l}{2}}=\frac{1}{2}\big[\sum\limits_{i=0}^{[\frac{n}{2}]}
e_{2i}(1-h)^{i}+\sqrt{1-h}\sum\limits_{j=0}^{[\frac{n-1}{2}]} e_{2j+1}(1-h)^{j}\big]\equiv C_{1}+C_{2},\end{equation} where
\begin{equation}\label{eq:4.14}
 \begin{array}{ll}
 C_{1}=\frac{1}{2}\sum\limits_{i=0}^{[\frac{n}{2}]}e_{2i}(1-h)^{i}
 =-\frac{1}{2}\sum\limits_{i=0}^{[\frac{n}{2}]}a^{-}_{2i}\Gamma_{2i+1,0}(1-h)^{i}&
        \hbox{}\\\\
    ~~~=L(a_{0}^{-},a_{2}^{-},\cdots,a_{2[\frac{n}{2}]}^{-})
          +L(a_{2}^{-},a_{4}^{-},\cdots,a_{2[\frac{n}{2}]}^{-})h+\cdots-
          \frac{1}{2}a^{-}_{2[\frac{n}{2}]}\Gamma_{2[\frac{n}{2}]+1,0}(-1)^{[\frac{n}{2}]}h^{[\frac{n}{2}]},
              \end{array}\end{equation}
 \begin{equation}\label{eq:4.15} \begin{array}{ll}
 2C_{2}=\sqrt{1-h}\sum\limits_{j=0}^{[\frac{n-1}{2}]}e_{2j+1}(1-h)^{j}
 =\sqrt{1-h}\sum\limits_{j=0}^{[\frac{n-1}{2}]}a^{-}_{2j+1}\Gamma_{2j+2,0}(1-h)^{j}&
        \hbox{}\\\\
   ~~~~~=\sqrt{1-h}\sum\limits_{j=0}^{[\frac{n-1}{2}]}a^{-}_{2j+1}\Gamma_{2j+2,0}
       \sum\limits_{k=0}^{j}C_{j}^{k}(-1)^{k}h^{k}&
        \hbox{}\\\\
  ~~~~~=\sqrt{1-h}\sum\limits_{k=0}^{[\frac{n-1}{2}]}\Big(\sum\limits_{j=k}^{[\frac{n-1}{2}]}a^{-}_{2j+1}\Gamma_{2j+2,0}
       C_{j}^{k}(-1)^{k}\Big)h^{k}. \end{array} \end{equation}
Let
\begin{equation}\label{eq:4.160}\sqrt{1-h}=\sum\limits_{i\geq0}\mu_{i}h^{i},~0<h\ll
1,\end{equation} where
\begin{equation}\label{eq:4.16}
\mu_{0}=1,~\mu_{1}=-\frac{1}{2},~\mu_{i}=-\frac{(2i-3)!!}{2^{i}},~i\geq2.
 \end{equation}
Denote
\begin{equation}\label{eq:4.17}
D_{k}=\sum\limits_{m=k}^{[\frac{n-1}{2}]}a^{-}_{2m+1}\Gamma_{2m+2,0}
       C_{m}^{k}(-1)^{k},~0\leq k\leq [\frac{n-1}{2}].
 \end{equation}
Then, by (\ref{eq:4.15}), (\ref{eq:4.160}) and (\ref{eq:4.17}) we have
\begin{equation}\label{eq:4.19}
2C_{2}=\sum\limits_{i=0}^{[\frac{n-1}{2}]}(\sum\limits_{j=0}^{i}\mu_{i-j}D_{j})h^{i} +\sum\limits_{k\geq
1}\big(\sum\limits_{s=0}^{[\frac{n-1}{2}]}\mu_{k+s}D_{[\frac{n-1}{2}]-s}\big)h^{[\frac{n-1}{2}]+k}, ~0<h\ll 1.
 \end{equation}
From (\ref{eq:4.9}), (\ref{eq:4.14}), (\ref{eq:4.17}) and (\ref{eq:4.19}), we get
\begin{equation}\label{eq:4.20}
\frac{M^{-}(\frac{h}{2})}{-2\sqrt{1-h}}=\sum\limits_{i\geq0}\omega_{i}h^{i},~ 0<h\ll 1,\end{equation} where

\begin{equation}\label{eq:4.200}  \begin{array}{ll}
      \omega_{i}=L(a^{-}_{0},a^{-}_{1},\cdots,a^{-}_{n}), ~ 0\leq i\leq [\frac{n-1}{2}],&
        \hbox{}\\\\
    \omega_{[\frac{n-1}{2}]+1}=\left\{
        \begin{array}{ll}
          \frac{1}{2}\sum\limits_{s=0}^{\frac{n-1}{2}}\mu_{1+s}D_{\frac{n-1}{2}-s},~n ~{\rm odd},&
           \hbox{}\\\\
            \frac{1}{2}\sum\limits_{s=0}^{\frac{n-2}{2}}\mu_{1+s}D_{\frac{n-2}{2}-s}
            +\frac{(-1)^{\frac{n}{2}+1}}{2}\Gamma_{n+1,0}a_{n}^{-},~n~{\rm even},& \hbox{}
        \end{array}
      \right.&
        \hbox{}\\\\
    \omega_{[\frac{n-1}{2}]+k}=\frac{1}{2}\sum\limits_{s=0}^{[\frac{n-1}{2}]}\mu_{k+s}D_{[\frac{n-1}{2}]-s},~k\geq 2.
                       \end{array}\end{equation}
Hence, it follows from (\ref{eq:3.6}) and (\ref{eq:4.20}) that
\begin{equation}\label{eq:4.21}
\frac {M(\frac{h}{2})}{-2\sqrt{1-h}}=\sum\limits_{i\geq0}(\widetilde{v}_{i}h^{i}+v_{i}^{*}h^{i+1}{\rm ln}h),~ 0<h\ll 1,\end{equation} where by (\ref{eq:3.611}) and (\ref{eq:4.200})

\begin{equation}\label{eq:4.210}  \begin{array}{ll}
      \widetilde{v}_{0}=a^{+}_{0}+L(a_{2}^{+},\cdots,a^{+}_{2[\frac{n}{2}]},p_{00}^{+},\cdots,p^{+}_{0,2[\frac{n-1}{2}]})
    +L(a_{0}^{-},a_{1}^{-},\cdots,a^{-}_{n}),&
        \hbox{}\\\\
    v^{*}_{0}=\frac{1}{4}p^{+}_{00},&
        \hbox{}\\\\
    \widetilde{v}_{1}=-\frac{a^{+}_{2}}{3}+L(a_{4}^{+},\cdots,a^{+}_{2[\frac{n}{2}]},p_{00}^{+},
    \cdots,p^{+}_{0,2[\frac{n-1}{2}]})+L(a_{0}^{-},a_{1}^{-},\cdots,a^{-}_{n}),&
        \hbox{}\\\\
    v^{*}_{1}=\frac{1}{4}p^{+}_{02}\beta_{1}^{*}+L(p^{+}_{00}),&
        \hbox{}\\\\
    \displaystyle \cdots\cdots &
    \hbox{}\\\\
   \widetilde{v}_{[\frac{n-1}{2}]}=\left\{
        \begin{array}{ll}
          \frac{a^{+}_{n-1}}{n}(-1)^{\frac{n-1}{2}}+L(p_{00}^{+},p_{02}^{+},\cdots,p^{+}_{0,n-1})
            +L(a_{0}^{-},a_{1}^{-},\cdots,a^{-}_{n}),~n~{\rm odd},&
           \hbox{}\\\\
             \frac{a^{+}_{n-2}}{n-1}(-1)^{\frac{n}{2}-1}+L(a^{+}_{n},p_{00}^{+},p_{02}^{+},\cdots,p^{+}_{0,n-2})
            +L(a_{0}^{-},a_{1}^{-},\cdots,a^{-}_{n}),~n~{\rm even},& \hbox{}
        \end{array}
      \right.&
        \hbox{}\\\\\\
     v^{*}_{[\frac{n-1}{2}]}=\frac{1}{4}p^{+}_{0,2[\frac{n-1}{2}]}\beta_{[\frac{n-1}{2}]}^{*}
       +L(p^{+}_{00},p_{02}^{+},\cdots,p^{+}_{0,2[\frac{n-1}{2}]-2}),&
    \hbox{}\\\\\\
  \widetilde{v}_{[\frac{n-1}{2}]+1}=\left\{
        \begin{array}{ll}
           \frac{1}{2}\sum\limits_{s=0}^{\frac{n-1}{2}}\mu_{1+s}D_{\frac{n-1}{2}-s}
            +L(p_{00}^{+},p_{02}^{+},\cdots,p^{+}_{0,n-1}),~n~{\rm odd},&
           \hbox{}\\\\
             \frac{a^{+}_{n}}{n+1}(-1)^{\frac{n}{2}}+L(p_{00}^{+},p_{02}^{+},\cdots,p^{+}_{0,n-2})
            +L(a_{n}^{-},a_{1}^{-},a_{3}^{-},\cdots,a^{-}_{n-1}),~n~{\rm even},& \hbox{}
        \end{array}
      \right.&
    \hbox{}\\\\
     v^{*}_{[\frac{n-1}{2}]+i}=O(|v^{*}_{0},v^{*}_{1},\cdots,v^{*}_{[\frac{n-1}{2}]}|),~i\geq 1,
     ({\rm obtained ~by~} (\ref{eq:3.6}) ~{\rm~and~}(\ref{eq:3.611}) ), &
    \hbox{}\\\\
     \widetilde{v}_{[\frac{n-1}{2}]+j}= \frac{1}{2}\sum\limits_{s=0}^{[\frac{n-1}{2}]}\mu_{j+s}D_{[\frac{n-1}{2}]-s}
            +L(p_{00}^{+},p_{02}^{+},\cdots,p^{+}_{0,2[\frac{n-1}{2}]}),~j\geq 2.          \end{array}\end{equation}

If $n$ is odd, by the formulas of $\widetilde{v}_{[\frac{n-1}{2}]+j},~j\geq 1$ in the above and $D_{k},~0\leq k\leq [\frac{n-1}{2}]$ in (\ref{eq:4.17}) we have
 $$\begin{array}{ll}
 \widetilde{v}_{\frac{n-1}{2}+j}=\frac{1}{2}\sum\limits_{s=0}^{\frac{n-1}{2}}\mu_{j+s}D_{\frac{n-1}{2}-s}+
 L(p_{00}^{+},p_{02}^{+},\cdots,p^{+}_{0,n-1})&
        \hbox{}\\\\
    ~~~~~~~~~=\frac{1}{2}\sum\limits_{s=0}^{\frac{n-1}{2}}\sum\limits_{m=\frac{n-1}{2}-s}^{\frac{n-1}{2}}
    \mu_{j+s}a^{-}_{2m+1}\Gamma_{2m+2,0}C_{m}^{\frac{n-1}{2}-s}(-1)^{\frac{n-1}{2}-s}
     +L(p_{00}^{+},p_{02}^{+},\cdots,p^{+}_{0,n-1})&
        \hbox{}\\\\\end{array}$$

   \begin{equation}\label{eq:4.22}~~~~~ \begin{array}{ll}
    ~~~~~~~~~=\frac{1}{2}\sum\limits_{m=0}^{\frac{n-1}{2}}\Big(\sum\limits_{s=\frac{n-1}{2}-m}^{\frac{n-1}{2}}
    \mu_{j+s}\Gamma_{2m+2,0}C_{m}^{\frac{n-1}{2}-s}(-1)^{\frac{n-1}{2}-s}\Big)a^{-}_{2m+1}
    +L(p_{00}^{+},p_{02}^{+},\cdots,p^{+}_{0,n-1})&
        \hbox{}\\\\ ~~~~~~~~~=\frac{1}{2}\sum\limits_{m=0}^{\frac{n-1}{2}}\Big(\sum\limits_{t=0}^{m}
    \mu_{j+\frac{n-1}{2}-t}\Gamma_{2m+2,0}C_{m}^{t}(-1)^{t}\Big)a^{-}_{2m+1}
    +L(p_{00}^{+},p_{02}^{+},\cdots,p^{+}_{0,n-1}),
            \end{array} \end{equation}
where $j\geq 1.$ For $1\leq j\leq\frac{n+1}{2}$, it follows from (\ref{eq:4.22})
\begin{equation}\label{eq:4.23}
     \displaystyle A_{1}\equiv\frac{\partial(\widetilde{v}_{\frac{n-1}{2}+1},\widetilde{v}_{\frac{n-1}{2}+2},
     \cdots,\widetilde{v}_{n})}
{\partial(a_{1}^{-},a^{-}_{3},\cdots,a^{-}_{n-2},a^{-}_{n})}
\end{equation}
$$=\left(
   \begin{array}{cccc}
     \frac{1}{2}\mu_{1+\frac{n-1}{2}}\Gamma_{20} & \frac{1}{2}\sum\limits_{t=0}^{1}\mu_{1+\frac{n-1}{2}-t}\Gamma_{40}C_{1}^{t}(-1)^t & \cdots  & \frac{1}{2}\sum\limits_{t=0}^{\frac{n-1}{2}}\mu_{1+\frac{n-1}{2}-t}\Gamma_{n+1,0}C_{\frac{n-1}{2}}^{t}(-1)^t\\
     \frac{1}{2}\mu_{2+\frac{n-1}{2}}\Gamma_{20} &  \frac{1}{2}\sum\limits_{t=0}^{1}\mu_{2+\frac{n-1}{2}-t}\Gamma_{40}C_{1}^{t}(-1)^t & \cdots & \frac{1}{2}\sum\limits_{t=0}^{\frac{n-1}{2}}\mu_{2+\frac{n-1}{2}-t}\Gamma_{n+1,0}C_{\frac{n-1}{2}}^{t}(-1)^t\\
     \cdots  & \cdots & \cdots & \cdots \\
     \frac{1}{2}\mu_{n}\Gamma_{20} & \frac{1}{2}\sum\limits_{t=0}^{1}\mu_{n-t}\Gamma_{40}C_{1}^{t}(-1)^t & \cdots & \frac{1}{2}\sum\limits_{t=0}^{\frac{n-1}{2}}\mu_{n-t}\Gamma_{n+1,0}C_{\frac{n-1}{2}}^{t}(-1)^t\\
            \end{array}
 \right).
 $$
 Let $c_i$ denote the $i$th column of $A_1,$ $a\times c_i$ (or $c_i/a$) indicate multiplying (or dividing) each element of the $i$th column by $a$, and $c_i+b\times c_j$ indicate that each element of the $j$th column times $b$ adds to the according element of the $i$th column, $a$ and $b$ are constants. We make elementary transformations to $A_1$ with the following steps
\par $S_1$. $c_i/(\frac{1}{2}\Gamma_{2i,0}),$ $i=1,2,\cdots,\frac{n+1}{2}.$
\par $S_2$. $c_i+(-1)\times c_1,$ $i=2,3,\cdots,\frac{n+1}{2}.$
\par $S_3$. $c_i+(-C_{i-1}^{1})\times c_2,$ $i=3,4,\cdots,\frac{n+1}{2};$ $(-1)\times c_2.$
\par $\cdots\cdots$
\par $S_\frac{n+1}{2}$. $c_{\frac{n+1}{2}}+(-\frac{n-1}{2})\times c_{\frac{n-1}{2}},$ $c_{\frac{n-1}{2}}/(-1)^{\frac{n-3}{2}}$ and $c_{\frac{n+1}{2}}/(-1)^{\frac{n-1}{2}}.$
\par\noindent Then $A_1$ becomes
\begin{equation}\label{eq:324}\tilde{A}_1\equiv\left(
   \begin{array}{ccccc}
     \mu_{\frac{n-1}{2}+1} & \mu_{\frac{n-1}{2}} & \cdots & \mu_{2} & \mu_{1}\\
     \mu_{\frac{n-1}{2}+2} & \mu_{\frac{n-1}{2}+1} & \cdots & \mu_{3} & \mu_{2}\\
     \cdots & \cdots & \cdots & \cdots & \cdots \\
     \mu_{n} & \mu_{n-1} & \cdots & \mu_{\frac{n-1}{2}+2} & \mu_{\frac{n-1}{2}+1}\\
            \end{array}\right).
 \end{equation}
By (\ref{eq:4.16}), we can prove Rank$(\tilde{A}_{1})=\frac{n+1}{2}$. See Appendix for the proof. Therefore, by using (\ref{eq:4.210}) and (\ref{eq:324}) we obtain
$$ \frac{\partial(\widetilde{v}_{0},\widetilde{v}_{1},\cdots,
   \widetilde{v}_{\frac{n-1}{2}},v^{*}_{0},v^{*}_{1},\cdots,v^{*}_{\frac{n-1}{2}},
   \widetilde{v}_{\frac{n-1}{2}+1},\cdots,\widetilde{v}_{n})}
{\partial(a^{+}_{0},a^{+}_{2},\cdots,a^{+}_{n-1},p^{+}_{00},p^{+}_{02},\cdots,p^{+}_{0,n-1}, a^{-}_{1},a^{-}_{3},\cdots,a^{-}_{n},a^{-}_{0},a^{-}_{2},\cdots,a^{-}_{n-1})}=$$

{\small $$  \left(
   \begin{array}{cccccccccccccccc}
     1 & * & \cdots & * & * & * & \cdots & * & * & * & \cdots & * & * & * & \cdots & *\\
     0 & -\frac{1}{3} & \cdots & * & * & * & \cdots & * & * & * & \cdots & * & * & * & \cdots & *\\
     \cdots & \cdots & \cdots & \cdots & \cdots & \cdots & \cdots & \cdots & \cdots & \cdots
     & \cdots & \cdots & \cdots & \cdots & \cdots & \cdots \\
     0 & 0 & \cdots & \frac{(-1)^{\frac{n-1}{2}}}{n} & * & * & \cdots & * & * & * & \cdots & * & * & * & \cdots & *\\
     0 & 0 & \cdots & 0 & \frac{1}{4} & 0 & \cdots & 0 & 0 & 0 & \cdots & 0 & 0 & 0 & \cdots & 0\\
     0 & 0 & \cdots & 0 & * & \frac{\beta_{1}^{*}}{4} & \cdots & 0 & 0 & 0 & \cdots & 0 & 0 & 0 & \cdots & 0\\
     \cdots & \cdots & \cdots & \cdots & \cdots & \cdots & \cdots & \cdots & \cdots & \cdots
     & \cdots & \cdots & \cdots & \cdots & \cdots & \cdots \\
     0 & 0 & \cdots & 0 & * & * & \cdots & \frac{\beta_{\frac{n-1}{2}}^{*}}{4} & 0 & 0 & \cdots & 0 & 0 & 0 & \cdots & 0\\
     0 & 0 & \cdots & 0 & * & * & \cdots & * &  &  &  &  & 0 & 0 & \cdots & 0\\
     \cdots & \cdots & \cdots & \cdots & \cdots & \cdots & \cdots & \cdots &  &
     A_{1} &  &  & \cdots & \cdots & \cdots & \cdots \\
     0 & 0 & \cdots & 0 & * & * & \cdots & * &  &  &  &  & 0 & 0 & \cdots & 0\\
            \end{array}
 \right)
 $$}

It is easy to see that the rank of this matrix is $n+\frac{n-1}{2}+2.$ Thus, we can choose $\widetilde{v}_{0},v^{*}_{0},\widetilde{v}_{1},v^{*}_{1},\cdots,
   \widetilde{v}_{\frac{n-1}{2}},$ $v^{*}_{\frac{n-1}{2}},\widetilde{v}_{\frac{n-1}{2}+1}\cdots,\widetilde{v}_{n}$ as free
   parameters such that
$$0<\widetilde{v}_{0}\ll v^{*}_{0}\ll \widetilde{v}_{1}\ll v^{*}_{1}\ll\cdots\ll
   \widetilde{v}_{\frac{n-1}{2}}\ll v^{*}_{\frac{n-1}{2}}$$ $$\ll \widetilde{v}_{\frac{n-1}{2}+1}\ll
   -\widetilde{v}_{\frac{n-1}{2}+2}\ll\cdots\ll (-1)^{\frac{n-1}{2}}\widetilde{v}_{n}\ll 1,$$
   or
$$0<-\widetilde{v}_{0}\ll -v^{*}_{0}\ll -\widetilde{v}_{1}\ll -v^{*}_{1}\ll\cdots\ll
   -\widetilde{v}_{\frac{n-1}{2}}\ll -v^{*}_{\frac{n-1}{2}}$$ $$\ll -\widetilde{v}_{\frac{n-1}{2}+1}\ll
  \widetilde{v}_{\frac{n-1}{2}+2}\ll\cdots\ll (-1)^{\frac{n+1}{2}}\widetilde{v}_{n}\ll 1.$$
Note that in (\ref{eq:4.210})
$$ v^{*}_{\frac{n-1}{2}+i}=O(|v^{*}_{0},v^{*}_{1},\cdots,v^{*}_{\frac{n-1}{2}}|),~i\geq 1.$$
Then, by (\ref{eq:4.21}) $\bar{M}(h)$ has $n+\frac{n+1}{2}$ positive simple zeros near $h=0.$ Hence, system (\ref{eq:100.201}) can have $n+\frac{n+1}{2}$ $(=n+[\frac{n+1}{2}]$ with $n$
odd) limit cycles near the homoclinic loop $L_{0}$.

If $n$ is even, by using the same method as above, we get the rank of the following matrix
$$\displaystyle \frac{\partial(\widetilde{v}_{0},\widetilde{v}_{1},\cdots,
   \widetilde{v}_{\frac{n-2}{2}},\widetilde{v}_{\frac{n}{2}},v^{*}_{0},v^{*}_{1},\cdots,v^{*}_{\frac{n-2}{2}},
   \widetilde{v}_{\frac{n}{2}+1},\widetilde{v}_{\frac{n}{2}+2},\cdots,\widetilde{v}_{n})}
{\partial(a^{+}_{0},a^{+}_{2},\cdots,a^{+}_{n},p^{+}_{00},p^{+}_{02},\cdots,p^{+}_{0,n-2}, a^{-}_{1},a^{-}_{3},\cdots,a^{-}_{n-1},a^{-}_{0},a^{-}_{2},\cdots,a^{-}_{n})}$$ is
$n+1+\frac{n}{2}$, which means this matrix has a full rank in row. Thus, we can choose
$$0<\widetilde{v}_{0}\ll v^{*}_{0}\ll \widetilde{v}_{1}\ll v^{*}_{1}\ll\cdots\ll
   \widetilde{v}_{\frac{n-2}{2}}\ll v^{*}_{\frac{n-2}{2}}$$ $$\ll \widetilde{v}_{\frac{n}{2}}\ll
   -\widetilde{v}_{\frac{n}{2}+1}\ll\cdots\ll (-1)^{\frac{n}{2}}\widetilde{v}_{n}\ll 1,$$
   or
$$0<-\widetilde{v}_{0}\ll -v^{*}_{0}\ll -\widetilde{v}_{1}\ll -v^{*}_{1}\ll\cdots\ll
   -\widetilde{v}_{\frac{n-2}{2}}\ll -v^{*}_{\frac{n-2}{2}}$$ $$\ll -\widetilde{v}_{\frac{n}{2}}\ll
  \widetilde{v}_{\frac{n}{2}+1}\ll\cdots\ll (-1)^{\frac{n}{2}+1}\widetilde{v}_{n}\ll 1.$$
Also by (\ref{eq:4.21}) and $ v^{*}_{\frac{n-2}{2}+i}=O(|v^{*}_{0},v^{*}_{1},\cdots,v^{*}_{\frac{n-2}{2}}|),~i\geq 1,$  the function $\bar{M}(h)$ has $n+\frac{n}{2}$
$(=n+[\frac{n+1}{2}]$ with $n$ even) positive simple zeros satisfying $0<h\ll 1.$ Therefore, system (\ref{eq:100.2}) can also have $n+[\frac{n+1}{2}]$ limit cycles near the homoclinic
loop $L_{0}$ in this case. That is, $N_{homoc}(n)\geq n+[\frac{n+1}{2}]$. The proof ends.

\par {\bf Proof of Theorem 1.2.} By Theorem 1.1, it suffices to prove that $Z(n)\leq n+[\frac{n+1}{2}]$ for $n=1,2,3$ and 4.
\par Let $n=1$. Then by Lemma 2.1 and (\ref{eq:224})
$$\bar M(h)=A\sqrt{1-h}+B(1-h)+Ch\varphi_0\Big(\sqrt{\frac{1-h}{h}}\Big).$$
where $A,B,C$ are constants. We have
$$\bar M(h)=\sqrt{1-h}[A+B\lambda+C\varphi_1(\lambda)],$$
where $\lambda=\sqrt{1-h}\in (0,1)$, and $\varphi_1(\lambda)=\frac{1-\lambda^2}{\lambda}\varphi_0\Big(\frac{\lambda}{\sqrt{1-\lambda^2}}\Big).$ Note that
$$\Big(\varphi_0\Big(\frac{\lambda}{\sqrt{1-\lambda^2}}\Big)\Big)'=
\frac{1}{(1-\lambda^2)^2}.$$ It is direct that
 $$[\varphi_1(\lambda)]'=-(1+\lambda^{-2})\varphi_0\Big(\frac{\lambda}{\sqrt{1-\lambda^2}}\Big)+\frac{1}{\lambda(1-\lambda^2)},$$
 and
 $$
[\varphi_1(\lambda)]''=\frac{2}{\lambda^3}[\varphi_0\Big(\frac{\lambda}{\sqrt{1-\lambda^2}}\Big)-\frac{\lambda}{1-\lambda^2}].$$ Note that $\frac{\lambda}{1-\lambda^2}=u\sqrt{1+u^2}$
if $u=\frac{\lambda}{\sqrt{1-\lambda^2}}$. We have $\varphi_0(u)<u\sqrt{1+u^2}$ since $(\varphi_0(u)-u\sqrt{1+u^2})'<0$ for $u>0.$ Therefore $[\varphi_1(\lambda)]''<0,$ which implies
that $\bar M(h)$ has at most two simple zeros in the interval $(0,1)$. That is, $Z(1)\leq 2$.

For the case of $n=2$, it follows from Lemma 2.1 and (\ref{eq:224}) that
$$\bar M(h)=\sqrt{1-h}[A_1+A_2\lambda+A_3\lambda^2+A_4\varphi_1(\lambda)],$$
where $\lambda=\sqrt{1-h}\in (0,1)$, $A_{1},A_{2},A_{3}$ and $A_{4}$ are constants. By simple computation, we get
$$\big[\varphi_1(\lambda)\big]'''=-\frac{6}{\lambda^4}\Big[\varphi_0\Big(\frac{\lambda}{\sqrt{1-\lambda^2}}\Big)
+\frac{\lambda(4\lambda^{2}-3)}{3(1-\lambda^2)^{2}}\Big]<0,$$ since
$$\Big[\varphi_0\Big(\frac{\lambda}{\sqrt{1-\lambda^2}}\Big)
+\frac{\lambda(4\lambda^{2}-3)}{3(1-\lambda^2)^{2}}\Big]'=\frac{1}{(1-\lambda^{2})^{2}}>0.$$ Hence, $\bar M(h)$ has at most three simple zeros in the interval (0,1). This means
$Z(2)\leq 3$.

For $n=3$, we have
$$\bar{M}(h)=\sqrt{1-h} \Big[A_0+A_1\lambda+A_2\lambda^2+A_3\lambda^3+(B_0+B_1\lambda^2)\varphi_1(\lambda)\Big]\equiv\sqrt{1-h}M_{1}(\lambda).$$
Hence,
$$\begin{array}{rl}M_{1}^{(4)}(\lambda)&=12B_1\varphi_1''(\lambda)+8B_1\lambda\varphi_1'''(\lambda)
+(B_0+B_1\lambda^2)\varphi_1^{(4)}(\lambda)\\\\
&=B_1\big(12\varphi_1''(\lambda)+8\lambda\varphi_1'''(\lambda)+\lambda^2\varphi_1^{(4)}(\lambda)\big)+B_0\varphi_1^{(4)}(\lambda).\end{array}$$ We can find
$$\varphi_1^{(4)}(\lambda)=\frac{24}{\lambda^5}\Big[\varphi_0\Big(\frac{\lambda}{\sqrt{1-\lambda^2}}\Big)
-\frac{\lambda(5\lambda^{4}-7\lambda^{2}+3)}{3(1-\lambda^2)^{3}}\Big],$$ and
$$12\varphi_1''(\lambda)+8\lambda\varphi_1'''(\lambda)+\lambda^2\varphi_1^{(4)}(\lambda)=-\frac{8}{(1-\lambda^2)^3}.$$ Thus,
$$M_{1}^{(4)}(\lambda)=\frac{1}{(1-\lambda^2)^3}[-8B_1+B_0(1-\lambda^2)^3\varphi_1^{(4)}(\lambda)]
\equiv \frac{1}{(1-\lambda^2)^3} M_{2}(\lambda).$$ Since
$$\big((1-\lambda^2)^3\varphi_1^{(4)}(\lambda)\big)'
=-\frac{24(1-\lambda^2)^2(\lambda^2+5)}{\lambda^6}\Big[\varphi_0\Big(\frac{\lambda}{\sqrt{1-\lambda^2}}\Big) +\frac{\lambda(17\lambda^{2}-15)}{3(1-\lambda^2)^2(\lambda^2+5)}\Big]
$$
and
\begin{equation}\label{eq:363}
     \Big[\varphi_0\Big(\frac{\lambda}{\sqrt{1-\lambda^2}}\Big)
+\frac{\lambda(17\lambda^{2}-15)}{3(1-\lambda^2)^2(\lambda^2+5)}\Big]' =\frac{16\lambda^{6}}{(1-\lambda^{2})^{3}(5+\lambda^{2})^{2}}>0,
\end{equation}
we have for $\lambda\in (0,1)$
$$M_{2}'(\lambda)\neq 0~{\rm if}~B_{0}\neq 0.$$
Therefore, $\bar M(h)$ has at most 5 zeros in the interval (0,1) for this case. It follows that $Z(3)\leq5.$

For $n=4$, by using the same method above, we have
\begin{equation}\label{eq:364}\begin{array}{rl}\bar{M}(h)&=\sqrt{1-h} \Big[A_0+A_1\lambda+A_2\lambda^2+A_3\lambda^3+A_4\lambda^4+(B_0+B_1\lambda^2)\varphi_1(\lambda)\Big]\\\\
&\equiv\sqrt{1-h}M_{3}(\lambda),\end{array}\end{equation} and
$$\begin{array}{rl}M_{3}^{(5)}(\lambda)&=20B_1\varphi_1'''(\lambda)+10B_1\lambda\varphi_1^{(4)}(\lambda)
+(B_0+B_1\lambda^2)\varphi_1^{(5)}(\lambda)\\\\
&=B_1\big(20\varphi_1'''(\lambda)+10\lambda\varphi_1^{(4)}(\lambda)+\lambda^2\varphi_1^{(5)}(\lambda)\big) +B_0\varphi_1^{(5)}(\lambda),\end{array}$$ where
$$\varphi_1^{(5)}(\lambda)=-\frac{120}{\lambda^6}\Big[\varphi_0\Big(\frac{\lambda}{\sqrt{1-\lambda^2}}\Big)
+\frac{\lambda(30\lambda^{6}-59\lambda^{4}+50\lambda^{2}-15)}{15(1-\lambda^2)^{4}}\Big],$$ and
$$20\varphi_1'''(\lambda)+10\lambda\varphi_1^{(4)}(\lambda)+\lambda^2\varphi_1^{(5)}(\lambda)
=-\frac{48\lambda}{(1-\lambda^2)^4}.$$ Thus,
\begin{equation}\label{eq:365}M_{3}^{(5)}(\lambda)=\frac{\lambda}{(1-\lambda^2)^4}\Big[-48B_1
+B_0\frac{(1-\lambda^2)^4}{\lambda}\varphi_1^{(5)}(\lambda)\Big] .\end{equation} Noting that
$$\Big[\frac{(1-\lambda^2)^4}{\lambda}\varphi_1^{(5)}(\lambda)\Big]'
=-\frac{960(7-\lambda^2)(1-\lambda^2)^2}{\lambda^9}M_{4}(\lambda),$$ where
$$M_{4}(\lambda)=\varphi_0\Big(\frac{\lambda}{\sqrt{1-\lambda^2}}\Big)
+\frac{\lambda(48\lambda^4-155\lambda^{2}+105)}{15(\lambda^2-7)(1-\lambda^2)^2},$$ we have
\begin{equation}\label{eq:366}\Big[\frac{(1-\lambda^2)^4}{\lambda}\varphi_1^{(5)}(\lambda)\Big]'<0\end{equation}
since
$$M_{4}(0)=0, \ \ M_{4}'(\lambda)=
\frac{16\lambda^{8}(1-\lambda^2)}{5(1-\lambda^2)^4(7-\lambda^2)^2}>0.$$ Therefore, by (\ref{eq:364}), (\ref{eq:365}) and (\ref{eq:366}), $\bar M(h)$ has at most 6 simple zeros in the
interval (0,1), which yields   $Z(4)\leq 6.$ The proof ends.

\par {\bf Proof of Theorem 1.3.} Following Theorem 1.1, we only need to prove $Z(n)\leq 2n+[\frac{n+1}{2}],~n\geq 5.$
We give the proof using the same method as in  [11]. Let $\lambda=\sqrt{1-h},~h\in (0,1).$ Then by (\ref{eq:224}) and Lemma 2.1, it follows that
$$\begin{array}{ll}\displaystyle\bar{M}(h)&=\sqrt{1-h} f_n(\sqrt{1-h})+g_{[\frac{n-1}{2}]}(1-h)I_{10}(h)\hbox{}\\\\
\displaystyle&=\lambda f_n(\lambda)+(1-\lambda^2)g_{[\frac{n-1}{2}]}(\lambda^2)\varphi_{0}\big(\frac{\lambda}{\sqrt{1-\lambda^2}}\big)\hbox{}\\\\
\displaystyle&=\lambda f_n(\lambda)+u(\lambda)\varphi_{0}\big(\frac{\lambda}{\sqrt{1-\lambda^2}}\big)\hbox{}\\\\
\displaystyle&\equiv M^{*}(\lambda), \end{array}$$ where $u(\lambda)=(1-\lambda^{2})u_{0},~u_{0}=g_{[\frac{n-1}{2}]}(\lambda^2).$ Hence,

$$ F\equiv \frac{d}{d\lambda}\big(\frac{M^*(\lambda)}{u}\big)=\Big(\frac{\lambda f_n(\lambda)}{u}+\varphi_{0}\big(\frac{\lambda}{\sqrt{1-\lambda^2}}\big)\Big)'
=\frac{\vartheta(\lambda)}{u_{0}^2(1-\lambda^{2})^{2}}, $$ with
\begin{equation}\label{eq:412}
\displaystyle\vartheta(\lambda)=(\lambda f_n)'u-\lambda f_nu'+u_{0}^2.\end{equation} By (\ref{eq:412}) and degree($u$)=$2[\frac{n+1}{2}]$, we find that the coefficient of the term of
$\vartheta(\lambda)$ of  degree $n+2[\frac{n+1}{2}]$ vanishes if $n$ is odd. Thus, degree($\vartheta$)$\leq n+2[\frac{n+1}{2}]-1$ if $n$ is odd, and degree($\vartheta$)$\leq
n+2[\frac{n+1}{2}]$ if $n$ is even. That is, degree($\vartheta$)$\leq 2n.$ Let $I=(0,1).$ Using notation $\sharp \{\lambda\in I|f(\lambda)=0\}$ to indicate the number of zeros of the
function $f$ in the interval $I$ taking into account their multiplicities,  we have $$\sharp \{\lambda\in I|u(\lambda)=0\}\leq [\frac{n-1}{2}],~~\sharp \{\lambda\in
I|\vartheta(\lambda)=0\}\leq 2n.$$ Therefore,
$$ \begin{array}{ll}
    \sharp \{\lambda\in I|M^{*}(\lambda)=0\}&\leq  \sharp \{\lambda\in
    I|u=0\}+ \sharp \{\lambda\in I|F=0,~u\neq 0\}+1
        \hbox{}\\\\
   \displaystyle&\leq \sharp \{\lambda\in
    I|u=0\}+ \sharp \{\lambda\in I|\vartheta(\lambda)=0\}+1
        \hbox{}\\\\
   &\leq\displaystyle 2n+[\frac{n+1}{2}].
                       \end{array}$$
The proof is completed.
$$~$$
\\\\\\
\noindent {\bf\Large Appendix}
\par {\bf The proof of Rank($\tilde{A}_1)=\frac{n+1}{2}$ (for $n$ odd)}. For simplicity, we take $n=7$ for example, since the method used here can be generalized for any $n\geq 1.$ Thus, by (\ref{eq:4.16}) and (\ref{eq:324}),
$$\tilde{A}_1=\left(
   \begin{array}{cccc}
     -5!!/2^4 & -3!!/2^3  & -1/2^2  & -1/2\\
     -7!!/2^5 & -5!!/2^4  & -3!!/2^3  & -1/2^2\\
     -9!!/2^6 & -7!!/2^5  & -5!!/2^4  & -3!!/2^3  \\
     -11!!/2^7  & -9!!/2^6  & -7!!/2^5  & -5!!/2^4
            \end{array}
 \right).
 $$
Let $r_i(c_i$) denote the $i$th row(column) of $\tilde{A}_1$ for $i=1,2,3,4.$ Now we making elementary transformations to $\tilde{A}_1$. Multiplying $r_i$ by $-2^i$ for $i=1,2,3,4,$
and $c_j$ by $2^{4-j}$ for $j=1,2,3,$ gives

$$\tilde{A}_1\sim \bar{A}_1\equiv\left(
   \begin{array}{cccc}
     5!! & 3!!  & 1  & 1\\
     7!! & 5!!  & 3!!  & 1\\
     9!! & 7!!  & 5!!  &3!! \\
     11!!  & 9!!  & 7!!  & 5!!
            \end{array}
 \right).
 $$
Further, for $\bar{A}_1$, by adding $r_i$ times $-(2i-1)$ to $r_{i+1}$ for $i=3,2,1,$ and dividing $c_j$ by $2(4-j)$ for $j=1,2,3,$ we have
$$\bar{A}_1\sim\left(
   \begin{array}{cccc}
     * & *  & *  & 1\\
     5!! & 3!!  & 1  & 0\\
     7!! & 5!!  & 3!!  & 0 \\
     9!!  & 7!!  & 5!!  & 0
            \end{array}
 \right).
 $$
Using the same method as above, we also have
$$\left(
   \begin{array}{cccc}
     * & *  & *  & 1\\
     5!! & 3!!  & 1  & 0\\
     7!! & 5!!  & 3!!  & 0 \\
     9!!  & 7!!  & 5!!  & 0
            \end{array}
 \right)\sim
 \left(
   \begin{array}{cccc}
     * & *  & *  & 1\\
     * & *  & 1  & 0\\
     5!! & 3!!  & 0  & 0 \\
     7!!  & 5!!  & 0 & 0
            \end{array}
 \right)\sim
 \left(
   \begin{array}{cccc}
     * & *  & *  & 1\\
     * & *  & 1  & 0\\
     * & 1  & 0  & 0 \\
     5!!  &0  & 0  & 0
            \end{array}
 \right) .$$
Therefore, Rank($\tilde{A}_1)=4.$ The proof ends.

\end{document}